\documentclass[a4paper,10pt]{article}
\usepackage{stmaryrd}
\usepackage{amsfonts}
\usepackage{bbm}
\usepackage{amscd}
\usepackage{mathrsfs}
\usepackage{latexsym,amssymb,amsmath,amscd,amscd,amsthm,amsxtra,xypic}
\usepackage[dvips]{graphicx}
\usepackage[all]{xy}
\usepackage{ulem}

\newtheorem{thm}{Theorem}[section]
\newtheorem{lem}[thm]{Lemma}
\newtheorem{cor}[thm]{Corollary}
\newtheorem{pro}[thm]{Proposition}
\newtheorem{ex}[thm]{Example}
\newtheorem{rmk}[thm]{Remark}
\newtheorem{defi}[thm]{Definition}

\newcommand{\ve}{\mathrm{v}}
\newcommand {\emptycomment}[1]{} 

\newcommand{\be }{\begin{eqnarray*}}
\newcommand{\ee }{\end{eqnarray*}}
\newcommand{\ad}{\mathrm{ad}}
\newcommand{\add}{\frka\frkd}

\newcommand{\pf}{\noindent{\bf Proof.}\ }

\newcommand{\huaK}{\mathcal{K}}

\newcommand{\frka}{\mathfrak a}

\newcommand{\frkd}{\mathfrak d}

\newcommand{\frkg}{\mathfrak g}
\newcommand{\frkh}{\mathfrak h}

\newcommand{\frkk}{\mathfrak k}

\newcommand{\g}{\mathfrak g}

\newcommand{\h}{\mathbbm h}

\newcommand{\frkD}{\mathfrak D}

\def\gpd{\,\lower1pt\hbox{$\longrightarrow$}\hskip-.24in\raise2pt
         \hbox{$\longrightarrow$}\,}

\def\qed{\hfill ~\vrule height6pt width6pt depth0pt}

\newcommand{\str}{\mathrm{string}}

\newcommand{\ie}{\mathbbm{i}}

\newcommand{\dM}{\mathrm{d}}

\newcommand{\Hom}{\mathrm{Hom}}
\newcommand{\Der}{\mathrm{Der}}

\newcommand{\inn}{\mathrm{inn}}
\newcommand{\SDer}{\mathrm{SDer}}
\newcommand{\SAut}{\mathrm{SAut}}

\newcommand{\Ad}{\mathrm{Ad}}
\newcommand{\Aut}{\mathrm{Aut}}

\newcommand{\Img}{\mathrm{Im}}

\newcommand{\End}{\mathrm{End}}
\newcommand{\Sinn}{\mathrm{Sinn}}
\newcommand{\SInn}{\mathrm{SInn}}
\newcommand{\Inn}{\mathrm{Inn}}
\newcommand{\DER}{\mathrm{DER}}
\newcommand{\Ha}{\mathrm{H}}

\begin{document}
\title{
 {Integration of  derivations for Lie $2$-algebras
 \thanks
 {
Research supported by NSFC (11101179,11471139) and NSF of Jilin Province (20140520054JH).
 }
} }
\author{Honglei Lang$^1$, Zhangju Liu$^1$, Yunhe Sheng$^2$\\
$^1$School of Mathematics and LMAM, Peking University,  Beijing
100871, China\\
$^2$School of Mathematics, Jilin University,
 Changchun 130012,  China\\\vspace{3mm} email: hllang@pku.edu.cn,~liuzj@pku.edu.cn,~ shengyh@jlu.edu.cn\\}
\date{}
\maketitle

\begin{abstract}
In this paper, for a Lie 2-algebra $\g$, we construct the automorphism 2-group $\Aut(\g)$, which turns out to be an integration of the derivation Lie 2-algebra $\Der(\g)$.
\end{abstract}

\section{Introduction}

 Recently, people have paid a lot of attention to higher categorical
structures with motivation from string theory. One way to provide
higher categorical structures is by categorifying existing
mathematical concepts. One of the simplest higher structures is a
$2$-vector space, which is a categorified  vector space. If we
further put Lie algebra structures on $2$-vector spaces, then we
obtain  Lie $2$-algebras \cite{baez:2algebras}. The
Jacobi identity is replaced by a natural transformation, called the
Jacobiator, which also satisfies some coherence laws of its own. It is well-known that the category of Lie 2-algebras is equivalent to the category of  2-term
$L_\infty$-algebras
\cite{baez:2algebras}. The concept of an $L_\infty$-algebra (sometimes
called a strongly homotopy (sh) Lie algebra)  was originally
introduced as a
model for ``Lie algebras that satisfy Jacobi identity up to all
higher homotopies''. The structure of a Lie 2-algebra appears in many areas such as higher symplectic geometry
\cite{baez:classicalstring},
string theory \cite{baez:string} and Courant algebroids. 

To study nonabelian extensions of Lie 2-algebras, the authors introduced the notion of a derivation of a Lie 2-algebra $\g=(\g_0\oplus\g_{-1},d,[\cdot,\cdot],l_{3})$ in \cite{Sheng1}, and showed that there is a strict Lie 2-algebra $\Der(\g)$, in which the degree $-1$  part is the set of degree $-1$ derivations $\Der^{-1}(\g)$, and the degree 0 part is the set of degree 0 derivations $\Der^0(\g)$. Furthermore, one can also construct a Lie 3-algebra $\DER(\g)$ (called the derivation Lie 3-algebra), in which the degree $-2$  part is $\g_{-1}$, the degree $-1$  part is $\Der^{-1}(\g)\oplus\g_0$, and the degree 0 part is $\Der^0(\g)$. A nonabelian extension of a Lie 2-algebra $\g$ by a Lie 2-algebra $\frkh$ can be given by a homomorphism from $\g$ to $\DER(\frkh)$.

In this paper, we study the integration of the strict Lie 2-algebra $\Der(\g)$. Here, ``integration''
is meant in the sense in which a Lie algebra is integrated to a
corresponding Lie group. In recent years, there are many works about the integration
of  Lie-algebra-like structures, such as that of Lie algebroids
\cite{crainic, tz1}, of $L_\infty$-algebras \cite{getzler, henriques} and
of Courant algebroids \cite{LS,MT}. In particular, semidirect product Lie 2-algebras are integrated to strict Lie 2-groups via equivalence in \cite{Sheng2}, and a homomorphism between strict Lie 2-algebras is  integrated to a homomorphism between Lie 2-groups in \cite{Sheng3}. In the current paper, we construct a strict Lie 2-group $\Aut(\g)$ associated to automorphisms of $\g$, and show that $\Aut(\g)$ is an integration of the strict Lie 2-algebra $\Der(\g)$. This is a generalization of the classical fact that the automorphism Lie group of a Lie algebra is an integration of the Lie algebra of derivations. 
This extend recent work in two ways. For example,
if $\g$ is a strict Lie 2-algebra, with the crossed module of Lie groups $\mathbb G$ as its integration, we can construct a sub-Lie $2$-group $\SAut(\g)$ of $\Aut(\g)$. This turns out to be the differential of the automorphism Lie 2-group $\Aut(\mathbb{G})$ given in \cite{Norrie}; In case  $\g$ is abelian, i.e. both $[\cdot,\cdot]$ and $l_3$ are zero, $\End(\g)$ is a sub-Lie 2-algebra of $\Der(\g)$, where $\End(\g)$ is a generalization of the general linear Lie algebra from a vector space to a 2-term complex of vector spaces. In this case, our result also gives the integration of $\End(\g)$, which was already given in \cite{Sheng2}. See \cite{Urs} for more details about automorphisms of strict 2-groups.

The paper is organized as follows. In Section 2, we recall some basic
definitions regarding Lie 2-algebras, crossed modules of Lie algebras (Lie groups).  In Section 3, we work out the derivation Lie 2-algebra for the string Lie 2-algebra. We also  compare our derivations with homotopy derivations introduced in \cite{Lada} recently.  In Section 4, we construct a crossed module of Lie groups $(\Aut^{-1}(\g),\Aut^0(\g),\partial,\triangleright)$ associated to automorphisms of a Lie 2-algebra $\g$ (Theorem \ref{thm:2grp}). In Section 5,
we prove that the strict Lie 2-group $\Aut(\g)$ is an integration of the strict Lie 2-algebra $\Der(\g)$ (Theorem \ref{thm:main}). In addition, we give the expressions of exponential maps from $\Der^i(\g)$ to $\Aut^i(\g)$ for $i=0,-1$, which commute with the differentials.
In Section 6, we first consider a sub-Lie $2$-algebra $\SDer(\g)$ of $\Der(\g)$ and a sub-crossed module $\SAut(\g)$ of $\Aut(\g)$. Then we give the relation between $\SAut(\g)$ and $\Aut(\mathbb G)$ for a strict Lie 2-algebra $\g$, where $\mathbb G$ is the crossed module integrating $\g$. Finally, we explore the inner automorphism Lie 2-group of a Lie $2$-algebra.


\section*{Acknowledgement }
 We give our warmest thanks to referees for very helpful comments that improve the paper.

\section{Preliminaries}
$\bullet$ {\bf Lie 2-algebras and crossed modules of Lie algebras}

A   Lie 2-algebra\footnote{In this paper, all the Lie 2-algebras are semistrict. } is equivalent to a 2-term $L_\infty$-algebra \cite{baez:2algebras}.
The notion of $L_\infty$-algebras was introduced by Schlessinger and Stasheff in \cite{Stasheff}. See \cite{LadaMarkl} for more information.

\begin{defi}\label{defi:Lie 2}
A {\bf Lie $2$-algebra} is a graded vector space $\g= \g_{-1}
\oplus \g_{0}$, together with a differential  $d:\g_{-1}\longrightarrow\g_0$, a skew-symmetric bilinear map $[\cdot,\cdot]:\g_i\wedge\g_j\longrightarrow\g_{i+j},~-1\leq i+j\leq0$ and a skew-symmetric trilinear map
$l_3:\wedge^3\frkg_0\longrightarrow \g_{-1}$
satisfying the following equalities
\begin{enumerate}
\item[\rm(a)] $d [x,a]=[x,d a]$,\quad $[d a,b]=[a,d b]$,
\item[\rm(b)]$[[x,y],z]+c.p.=-d l_3(x,y,z)$,\quad $[[x,y],a]+c.p.=-l_3(x,y,d a)$,
\item[\rm(c)]$l_3([x,y],z,w)+c.p.=[l_3(x,y,z),w]+c.p.$,
\end{enumerate}
for all $x,y,z,w\in \g_0,~a,b\in \g_{-1}$, where  $c.p.$ means cyclic permutation. If $l_3=0$, $\g$ is called a
{\bf strict Lie $2$-algebra}.
\end{defi}
We denote a Lie $2$-algebra by $(\g_0\oplus \g_{-1},d,[\cdot,\cdot],l_{3})$, or simply by $\g$ if there is no confusion.


\begin{defi}\label{defi:Lie 2 homo}
Let $\g=(\g_0\oplus \g_{-1},d,[\cdot,\cdot],l_{3})$ and $\g{'}=(\g_0'\oplus\g_{-1}',d{'},[\cdot,\cdot]',l'_3)$ be Lie $2$-algebras. A {\bf Lie $2$-algebra homomorphism} $A:\g \longrightarrow{\g'}$ consists of
\begin{enumerate}
\item[$\bullet$] two linear maps $A_{0}:\g_{0}\longrightarrow{\g'_{0}}$ and $A_{1}:\g_{-1}\longrightarrow{\g'_{-1}}$, satisfying $d'\circ{A_1}=A_{0}\circ{d}$,
\item[$\bullet$] one bilinear map $A_{2}: \g_{0}\wedge{\g_{0}}\longrightarrow{\g'_{-1}}$,
\end{enumerate}
such that the following equalities hold for all $x,y,z\in{\g_{0}}, a\in{\g_{-1}}$ :
\begin{enumerate}
\item[\rm(i)] $A_{0}[x,y]-[A_{0}(x),A_{0}(y)]'=d'A_{2}(x,y)$,
\item[\rm(ii)] $A_{1}[x,a]-[A_{0}(x),A_{1}(a)]'=A_{2}(x,da)$,
\item[\rm(iii)] $[A_{0}(x),A_{2}(y,z)]'+c.p.+l'_{3}(A_{0}(x),A_{0}(y),A_{0}(z))=A_{2}([x,y],z)+c.p.+A_{1}(l_{3}(x,y,z))$.
\end{enumerate}
  It is called a {\bf strict
homomorphism} if $A_{2}=0$.
\end{defi}

\begin{defi}
A {\bf crossed module of Lie algebras} is a quadruple
$(\frkh_1,\frkh_0,\varphi,\phi)$, which we denote by $\h$, where
$\frkh_1$ and $\frkh_0$ are Lie algebras,
$\varphi:\frkh_1\longrightarrow\frkh_0$ is a Lie algebra homomorphism and
$\phi:\frkh_0\longrightarrow\Der(\frkh_1)$ is an action of $\frkh_0$ on $\frkh_1$ by derivations, such
that
$$
\varphi(\phi_x(a))=[x,\varphi(a)]_{\frkh_0},\quad \phi_{\varphi(a)}(b)=[a,b]_{\frkh_1}.
$$
\end{defi}

It is well-known that there is a one-to-one correspondence between strict Lie 2-algebras and
crossed modules of Lie algebras. In short, the formula for the correspondence can be given as follows: A strict Lie 2-algebra $\g_{-1}\stackrel{d}{\longrightarrow} \g_0$
gives rise to a crossed module with $\frkh_1=\g_{-1}$ and
$\frkh_0=\g_0$, where the Lie brackets are given by:
\begin{eqnarray*}
~[a,b]_{\frkh_1}=[d(a),b],\quad
~[x,y]_{\frkh_0}=[x,y],\quad\forall~x,y\in
\g_0,~a,b\in
\g_{-1}
\end{eqnarray*}
and $\varphi=d$, $\phi:\frkh_0\longrightarrow\Der(\frkh_1)$ is given by
$ \phi_x(a)=[x,a]. $ The strict Lie 2-algebra structure gives the Jacobi identities for
$[\cdot,\cdot]_{\frkh_1}$ and $[\cdot,\cdot]_{\frkh_0}$, and various other
conditions for crossed modules.

Conversely, a crossed module $(\frkh_1,\frkh_0,\varphi,\phi)$ gives rise
to a strict Lie 2-algebra with $d=\varphi$, $\g_{-1}=\frkh_1$ and $\g_0=\frkh_0$, and
$[\cdot,\cdot]$ is given by:
\begin{eqnarray*}
~[a,b]\triangleq0,\quad
~ [x,y]\triangleq[x,y]_{\frkh_0},\quad
~[x,a]=-[a,x]\triangleq\phi_x(a).
\end{eqnarray*}
$\bullet$ {\bf Lie 2-groups and crossed modules of Lie groups}

A group is a monoid where every element has an inverse. A 2-group is
a monoidal category where every object has a weak inverse and every
morphism has an inverse. Denote the category of smooth manifolds and
smooth maps by $\rm Diff$, a (semistrict) Lie 2-group is  a 2-group
in $\rm DiffCat$, where  $\rm DiffCat$ is the 2-category consisting
of categories, functors, and natural transformations in $\rm Diff$.
For more details, see \cite{baez:2gp}. Here we only give the
definition of strict Lie 2-groups.
\begin{defi}
A {\bf strict Lie $2$-group} is a Lie groupoid $C$
such that
\begin{itemize}
\item[\rm(1)]
The space of morphisms $C_1$ and the space of objects $C_0$ are Lie
groups.
\item[\rm(2)] The source and the target $s,t:C_1\longrightarrow C_0$, the identity assigning function $i:C_0\longrightarrow C_1$ and
 the composition $\circ:C_1\times_{C_0}C_1\longrightarrow C_1$ are all Lie group morphisms.
\end{itemize}
\end{defi}

It is known that strict Lie 2-groups can be described by
crossed modules of Lie groups.

\begin{defi}
A {\bf crossed module of Lie groups} is a quadruple $(H_1,H_0,\partial,\triangleright)$,
which we denote simply by $\mathbb H$, where $H_1$ and $H_0$ are Lie
groups, $\partial: H_1\longrightarrow H_0$ is a homomorphism, and
$\triangleright:H_0\longrightarrow \Aut(H_1)$ is an action of $H_0$ on
$H_1$ by automorphisms such that $\partial$ is $H_0$-equivariant:
\begin{equation}\label{cm g 1}
\partial(g\triangleright h)=g\partial(h)g^{-1},\quad \forall ~g\in H_0,~h\in H_1,
\end{equation}
and satisfies the so called Peiffer identity:
\begin{equation}\label{cm g 2}
\partial(h)\triangleright(h^\prime)=hh^\prime h^{-1},\quad\forall ~h,h^\prime\in
H_1.
\end{equation}
\end{defi}

It is well-known that there is a one-to-one correspondence between crossed modules of Lie
groups and strict Lie 2-groups. Roughly speaking, given a crossed
module $(H_1,H_0,\partial,\triangleright)$, there is a strict Lie
2-group for which $C_0=H_0$ and $C_1= H_0\times H_1$. In this strict Lie $2$-group,
the source and target maps $s,~t:C_1\longrightarrow C_0$ are given
by
$$
s(g,h)=g,\quad t(g,h)= \partial(h)g,
$$
the groupoid multiplication $ \cdot_\ve$ is given by
\begin{equation}\label{m v}
(g^\prime,h^\prime)\cdot_\ve(g,h) =(g, h^\prime h),\quad
\mbox{where} \quad g^\prime= \partial(h)g,
\end{equation}
the group multiplication $\cdot_\mathrm{h}$ is given by
\begin{equation}\label{m h}
(g,h)\cdot_\mathrm{h} (g^\prime,h^\prime)=(g
g^\prime,h(g\triangleright h^\prime)).
\end{equation}

\section{ Derivations of Lie 2-algebras}

Let $ \mathbb V:V_{-1}\stackrel{\dM}{\longrightarrow} V_0$ be a 2-term
complex of vector spaces. We can form a new 2-term complex of vector
spaces $\End(\mathbb V):\End^{-1}(\mathbb
V)\stackrel{\delta}{\longrightarrow} \End^0_\dM(\mathbb V)$ by
defining $\delta(\Theta)=\dM\circ \Theta+\Theta\circ\dM$ for any
$\Theta\in\End^{-1}(\mathbb V)$, where $\End^{-1}(\mathbb V)=\Hom(V_0,V_{-1})$ and
$$\End^0_\dM(\mathbb V)=\{X=(X_0,X_1)\in \End(V_0)\oplus \End(V_{-1})|~X_0\circ \dM=\dM\circ X_1\}.$$
There is a natural bracket operation $[\cdot,\cdot]_C$ on $\End(\mathbb V)$ given by the commutator:
\begin{eqnarray*}
  ~[(X_0,X_1),(Y_0,Y_1)]_C&=&(X_0\circ Y_0-Y_0\circ X_0,X_1\circ Y_1-Y_1\circ X_1), \\
  ~ [(X_0,X_1),\Theta]_C&=&X_1\circ \Theta-\Theta\circ X_0.
\end{eqnarray*}
Consequently, $(\End(\mathbb V),\delta,[\cdot,\cdot]_C)$ is a strict Lie 2-algebra. It plays the same role as $\End(V)$ for a vector space $V$. Its integration is given in \cite{Sheng2}.\vspace{3mm}

Let $\g=(\g_0\oplus\g_{-1},d,[\cdot,\cdot],l_3)$ be a Lie $2$-algebra.

\begin{defi}\label{defi:derivation}
A {\bf derivation\footnote{Strictly speaking, we should use the terminology ``a weak derivation'' since the derivation conditions hold up to homotopy. To be succinct, we just say ``a derivation''.} of
degree $0$} of $\mathfrak{g}$ is a triple $(X_0,X_1,l_{X})$, where $X=(X_0,X_1)\in\End^0_d(\g)$
and $l_{X}:\g_0 \wedge{\g_0}\rightarrow{\g_{-1}}$ is a linear map, such
that for all $x,y,z\in\g_0,a\in\g_{-1}$,
\begin{itemize}
\item[\rm(a)]  $dl_{X}(x,y)=X_0[x,y]-[X_0x,y]-[x,X_0y],$
\item[\rm(b)] $l_{X}(x,da)=X_1[x,a]-[X_0x,a]-[x,X_1a],$
\item[\rm(c)] $X_1l_{3}(x,y,z)=l_{X}(x,[y,z])+[x,l_{X}(y,z)]+l_{3}(X_0x,y,z)+c.p.(x,y,z).$
\end{itemize}
\end{defi}

Denote by $\Der^{0}(\g)$ the set of derivations of degree $0$ of
$\g$. Then we can obtain a $2$-vector space:
\begin{equation*}
\CD
   \Der(\g): \Der^{-1}(\g)\triangleq\End^{-1}(\g) @>\bar{d}>>\Der^{0}(\g),
\endCD
\end{equation*}
where $\bar{d}$ is given by
$\bar{d}(\Theta)=(\delta(\Theta),l_{\delta(\Theta)})$, in which
$\delta(\Theta)=(d\circ \Theta,\Theta\circ d)$ and
\[l_{\delta(\Theta)}(x,y)={\Theta}[x,y]-[x,\Theta{y}]-[{\Theta}x,y].\]
In addition, define
\begin{eqnarray}
\label{eq:bra01}\{(X,l_{X}),\Theta\}&=&[X,\Theta]_C,\\
\label{bra of der}
\{(X,l_X),(Y,l_Y)\}&=&([X,Y]_C,L_X(l_Y)-L_Y(l_X)),
\end{eqnarray}
where $[\cdot,\cdot]_C$ is the commutator bracket and for all $X=(X_0,X_1)\in\End(\g_0)\oplus\End(\g_{-1})$, $L_X:\Hom(\wedge ^k\g_0,\g_{-1})\longrightarrow\Hom(\wedge ^k\g_0,\g_{-1})$  is given by
$$L_X\omega(x_1,\cdots,x_k)=X_1\omega(x_1,\cdots,x_k)-\sum_{i=1}^k\omega (x_1,\cdots,X_0(x_i),\cdots,x_k).$$
\begin{thm}{\rm\cite{Sheng1}}\label{thm:Der(g)}
With the notations above, $(\Der(\g),\{\cdot,\cdot\})$ is a strict Lie
$2$-algebra. We call it the {\bf derivation Lie $2$-algebra} of $\g$.
\end{thm}

\begin{rmk}\label{rmk:der}{\rm


 Recently, the notion of a homotopy derivation was introduced in \cite{Lada} using the theory of operads.  See \cite[Proposition 4.1]{Lada} for precise formulas. Restricting to the 2-term case, for a Lie 2-algebra $\g=(\g_0\oplus\g_{-1},d,[\cdot,\cdot],l_3)$\footnote{In \cite{Lada}, the authors use an alternative definition of $L_\infty$-algebras. For their convention, $\g_0$ is of degree $-1$ and $\g_{-1}$ is of degree $-2$.}, we could obtain that a homotopy derivation of degree 0 is the same as the one given in Definition \ref{defi:derivation}, and a homotopy derivation of degree $-1$ is a $\Theta\in\End^{-1}(\g)$ such that
      \begin{eqnarray}
      \label{eq:r1}d\circ \Theta&=&0,\quad \Theta\circ d=0,\\
       \label{eq:r2}\Theta[x,y]&=&[\Theta(x),y]+[x,\Theta(y)].
      \end{eqnarray}
      There are also homotopy derivations of other degrees. However, since we are dealing with Lie 2-algebras, homotopy derivations of degree 0 and $-1$ are enough.

  It is straightforward to see that the Lie 2-algebra of homotopy derivations is a sub-Lie 2-algebra of the one given in Theorem \ref{thm:Der(g)}. In Subsection 6.1, we will consider strict derivations of strict Lie 2-algebras, which correspond to strict automorphisms of strict Lie 2-groups given in \cite{Norrie}. It turns out that for a strict Lie 2-algebra, homotopy derivations and strict derivations are not consistent. More explicitly, even for a strict Lie $2$-algebra, a homotopy derivation $(X,l_X)$ of degree $0$ can have a nonzero homotopy term ``$l_X$'', while the strict derivation of degree $0$ requires that $l_X=0$. Moreover, a homotopy derivation of degree $-1$ is more restrictive than a strict derivation of degree $-1$ since the former requires that both \eqref{eq:r1} and \eqref{eq:r2} hold, and  the latter only requires that \eqref{eq:r2} holds.
  Our derivation Lie 2-algebra can be viewed as a unification of them in the sense that the Lie $2$-algebra of homotopy derivations and the Lie $2$-algebra of strict derivations are both sub-Lie $2$-algebras of ours.}
  \end{rmk}

Let $(\Der(\g),\{\cdot,\cdot\})$ be the   derivation Lie $2$-algebra of $\g$, and denote the corresponding crossed module of Lie algebras by $(\Der^{-1}(\g),\Der^0(\g),\bar{d},\phi)$, where the action $\phi$ is given by \eqref{eq:bra01}, and
 the Lie bracket on $\Der^{-1}(\g)$  is given by
\begin{equation}\label{eq:bra of der0}
\{\Theta,\Theta'\}=\{\bar{d}(\Theta),\Theta'\}=[\delta(\Theta),\Theta']_C=\Theta\circ d\circ\Theta'-\Theta'\circ d\circ \Theta.
\end{equation}

Now we work out the derivation Lie 2-algebra of the string Lie 2-algebra $\str(\frkk)$. Let$(\mathfrak{k},[\cdot,\cdot]_\frkk)$ be a semisimple Lie algebra and $K$  the Killing form.  Consider the corresponding string Lie 2-algebra $\str(\frkk)=(\mathbb{R}\stackrel{0}\longrightarrow \mathfrak{k},[\cdot,\cdot],l_3)$. More precisely, $$[x,y]=[x,y]_\frkk,\  [x,r]=0,\
l_3(x,y,z)=K([x,y]_\frkk,z), \ \forall x,y,z\in \mathfrak{k},\ r\in \mathbb{R}.$$

\begin{pro} \label{pro:string}
For the string Lie $2$-algebra $\str(\frkk)$, $\Der^0(\str(\frkk))$ is isomorphic to the semidirect product Lie algebra  $\frkk\ltimes_{\ad^*}\frkk^*$. Furthermore, $ \Der^{-1}(\str(\frkk))=\frkk^* $, which is abelian, and the differential $\bar{d}$ is given by $\bar{d}(\Theta)=(0,0,-\mathfrak{D}\Theta)$, where $\mathfrak{D}$ is the coboundary operator on $\frkk$ with the coefficients in the trivial representation.
 \end{pro}
\proof
A straightforward computation shows that $(X_0,t,l_X)$, where $X_0\in\End(\frkk),t\in\mathbb R,l_X\in\wedge^2\frkk^*$, is a derivation of degree $0$ if and only if
$$
X_0\in\Der(\frkk),\quad \mathfrak{D} l_X(x,y,z)=tl_3(x,y,z)-l_3(X_0x,y,z)-l_3(x,X_0y,z)-l_3(x,y,X_0z).
$$
 Since $\frkk$ is semisimple, every derivation is an inner derivation, and thus skew-symmetric. Then we have
\begin{eqnarray*}
 l_3(X_0x,y,z)+c.p.&=&K([X_0x,y]_\frkk,z)+K([x,X_0y]_\frkk,z)+K([x,y]_\frkk,X_0z)\\
 &=&K([X_0x,y]_\frkk+[x,X_0y]_\frkk,z)-K(X_0[x,y]_\frkk,z)\\
 &=&0.
\end{eqnarray*}
Now the second condition reduces to $\mathfrak{D} l_X(x,y,z)=tl_3(x,y,z).$
However,  the Cartan 3-form $l_3$ could not be exact, which implies that $t=0$.  Therefore,  we have
$$
\Der^0(\str(\frkk))=\{(X_0,0,\omega)|X_0\in\Der(\frkk),~\omega\in\wedge^2\frkk^*~ \mbox{is a 2-cocycle}\}.
$$
Note that every $2$-cocycle is exact. For $(\ad_x,0,\frkD(\xi)),~(\ad_y,0,\frkD(\eta))\in \Der^0(\str(\frkk))$, we have
\begin{eqnarray*}
  \{(\ad_x,0,\frkD(\xi)), (\ad_y,0,\frkD(\eta))\}=(\ad_{[x,y]_\frkk},0,\frkD(\ad^*_x\eta-\ad^*_y\xi)),
\end{eqnarray*}
 which implies the map $(x,\xi)\longmapsto (\ad_x,0,\frkD(\xi))$ is an isomorphism between Lie algebras $\frkk\ltimes_{\ad^*}\frkk^*$ and $\Der^0(\str(\frkk))$. Other conclusions are obvious. The proof is finished. \qed \vspace{3mm}


More generally, we have

\begin{ex}\label{ex of der} {\rm
Let $(\frkk,[\cdot,\cdot]_\frkk)$ be a Lie algebra and $\phi:\frkk\longrightarrow \End(V)$ a representation of  $\frkk$ on $V$. Let $C^k(\frkk,V)$ be the set of $k$-cochains, i.e. $C^k(\frkk,V)=\{f:\wedge^k\frkk\longrightarrow V\}$. Given a Lie algebra
3-cocycle $l_3\in C^3(\frkk,V)$, we get a skeletal Lie $2$-algebra
$\g=(V\oplus {\frkk},d=0,[\cdot,\cdot],l_3)$, where $[\cdot,\cdot]$ is
defined by $$[x,y]=[x,y]_{\frkk},\ \ \ \ \
[x,u]=-[u,x]=\phi_x (u),\ \ \forall x,y\in \frkk,\ u\in V.$$

For $X=(X_0,X_1)\in \End(\frkk)\oplus\End(V)$,
and $l_X\in C^2(\frkk,V)$,  it is easy to check that
\begin{equation*}
(X,l_X)\in\Der^0(\g) \Longleftrightarrow \left\{\begin{array}{rcll}X&\in& \Der(\frkk\ltimes_\phi V),\\
\mathfrak{D} l_X&=&L_X(l_3),\end{array}\right.
\end{equation*}
where $\mathfrak{D}: C^k(\frkk,V)\rightarrow C^{k+1}(\frkk,V)$ is the
 coboundary operator on $\frkk$ with coefficients in $V$, and $\Der(\frkk\ltimes_\phi V)$ is the Lie
algebra of derivations of the semidirect product Lie algebra $\frkk\ltimes_\phi
V$.
Moreover, we
have $\Der^{-1}(\g)=C^1(\frkk,V)$ and the map $\bar{d}:
\Der^{-1}(\g)\rightarrow \Der^0(\g)$ is given by
$\bar{d}(\Theta)=(0,0,-\mathfrak{D}(\Theta))$.}
\end{ex}

A {\bf representation} of a Lie $2$-algebra $\g$ on a 2-term complex of vector spaces  $ \mathbb V:V_{-1}\stackrel{\dM}{\longrightarrow} V_0$ is a homomorphism from $\g$ to $\End(\mathbb V)$.
The {\bf adjoint representation} $\add=(\ad_0,\ad_1,\ad_2)$ of $\g$ on itself is given by
$$\ad_0(x)=[x,\cdot],\quad \ad_1(a)=[a,\cdot],\quad \ad_2(y,z)=-l_3(y,z,\cdot),\ \ \  \forall x,y,z\in \g_0,~a\in\g_{-1}.$$

For any Lie 2-algebra $\g$, the derivation Lie 2-algebra $\Der(\g)$ is strict. However, there is still a Lie 2-algebra homomorphism from $\g$ to $\Der(\g)$.
Define $$\overline{\add}_0:\g_0\longrightarrow\Der^0(\g),~\overline{\add}_1:\g_{-1}\longrightarrow\Der^{-1}(\g),~\overline{\add}_2:\wedge^2\g_0\longrightarrow\Der^{-1}(\g)$$   by
\begin{equation}\label{adjoint homo}
\overline{\add}_0(x)=(\ad_0(x),l_3(x,\cdot,\cdot)),\ \ \overline{\add}_1(a)=\ad_1(a),
\ \ \overline{\add}_2(y,z)=\ad_2(y,z).
\end{equation}
It is straightforward to obtain that
\begin{lem}
With the above notations, $\overline{\add}=(\overline{\add}_0,\overline{\add}_1,\overline{\add}_2)$ is a Lie $2$-algebra homomorphism from  $\g$ to $\Der(\g)$.
\end{lem}

\begin{rmk}{\rm

  In the classical case, for a semisimple Lie algebra  $\frkk$, $\ad:\frkk\longrightarrow\Der(\frkk)$ is an isomorphism of Lie algebras. In the case of Lie 2-algebras, we can not expect such a result.
    \begin{itemize}\item[$\bullet$]  Even for the simple string Lie 2-algebra $\str(\frkk)$ given in Proposition \ref{pro:string}, $\overline{\add}_0:\frkk\longrightarrow\Der^0(\str(\frkk))$ is not surjective.

    \item[$\bullet$] The image of $\overline{\add}_0$ is not closed in $\Der^0(\g)$. Explicitly, we have
\begin{eqnarray*}
\{\overline{\add}_0(x),\overline{\add}_0(y)\}&=&\overline{\add}_0([x,y])+\bar{d}l_3(x,y,\cdot).
\end{eqnarray*}
Thus, we can not define the inner derivation $\inn^0(\g)$ simply by $\Img \overline{\add}_0$.
 \end{itemize}}
\end{rmk}

\begin{defi}\label{inner derivation}
   Define the set of inner derivations of degree $0$, which is denoted by $\inn^0(\g)$, by
\begin{equation}\label{eq:inn0}
  \inn^0(\g)=\Img \overline{\add}_0+ \Img\bar{d}.
\end{equation}
\end{defi}
Then we can get a Lie $2$-algebra $\inn(\g)$ given by
\begin{equation}\label{eq:inner}
\CD
   \inn(\g): \inn^{-1}(\g)\triangleq \End^{-1}(\g)@>\bar{d}>>\inn^{0}(\g),
\endCD
\end{equation}
which  we call   the {\bf inner derivation Lie $2$-algebra} of $\g$.
\begin{rmk}\label{coho}{\rm In this remark, we give some reasons for our definition of $\inn^{0}(\g)$.
From the homological viewpoint, $\Der^0(\g)$ is the
set of $1$-cocycles of the Lie $2$-algebra $\g$ with respect to the
adjoint representation $\add$, and $\inn^0(\g)$ is the set of $1$-coboundaries. See \cite{Lang} for details. For a Lie algebra $\frkk$, $\inn(\frkk)=\Img(\ad)$. However, for a Lie 2-algebra $\g$, they are not the same.  Furthermore, it is straightforward to see that $\inn(\g)$ is  an ideal of  $\Der(\g)$.

}
\end{rmk}

\begin{rmk}\label{rmk:innerder}{\rm
 Parallel to Remark \ref{rmk:der}, we compare our inner derivations with homotopy inner derivations given in \cite{Lada}.
For a semistrict Lie 2-algebra, a homotopy inner derivation of degree 0 is given by $(\ad_0(x),l_3(x,\cdot,\cdot))$ for $x\in\g_0$, and a homotopy inner derivation of degree $-1$ is given by $\ad_1(a)$ for $a\in\g_{-1}$ satisfying $d(a)=0$. See \cite[Proposition 4.4]{Lada} for general formulas.

The Lie 2-algebra of homotopy inner derivations is a sub-Lie 2-algebra of $\inn(\g)$. For a strict Lie 2-algebra,  homotopy inner derivations of degree $-1$ are not consistent with strict inner derivations  due to the condition $d(a)=0$. Our inner derivation Lie 2-algebra $\inn(\g)$ can also be viewed as a unification of them. }
\end{rmk}

\emptycomment{is it necessary to define Lie 2-subalgebras and ideal. If we need to verify it is an ideal, please tell me. So we get an exact sequence of Lie $2$-algebras
$$0\rightarrow \Inn(\g)\stackrel{}\longrightarrow \Der(\g)\stackrel{}\longrightarrow (0\rightarrow \Ha^1(\g))\rightarrow 0,$$
where $\Ha^1(\g)$ is the first cohomology group of $\g$ with respect to the adjoint action.
}

\section{The automorphisms of a Lie 2-algebra}

In this section, we construct a Lie 2-group $\Aut(\g)$ associated to the automorphisms of a Lie $2$-algebra $\g=(\g_0\oplus\g_{-1},d,[\cdot,\cdot],l_3)$.

Clearly, $\ie=(I_{\g_0},I_{\g_{-1}},0):\g\longrightarrow\g$ is a Lie $2$-algebra homomorphism, which is called the identity homomorphism. Let $A:\frkg\to \frkg'$ and $B:\frkg'\to \frkg''$ be two Lie 2-algebra homomorphisms. Then their composition $B\diamond A:\frkg\to \frkg''$
is also a homomorphism defined as $(B\diamond A)_0=B_0\circ A_0:\frkg_0\to \frkg''_0$, $(B\diamond A)_1=B_1\circ A_1:\frkg_{-1}\to \frkg''_{-1}$
and
$$(B\diamond A)_2=B_2\circ (A_0\times A_0)+B_1\circ A_2:\frkg_0\wedge \frkg_0\to \frkg''_{-1}.$$
A homomorphism $A:\frkg\to \frkg'$ is called an {\bf isomorphism} if there exists a homomorphism $A^{-1}:\frkg'\to \frkg$ such that
the compositions $A^{-1}\diamond A:\frkg\to \frkg$ and $A\diamond  A^{-1}:\frkg'\to \frkg'$ are both identity homomorphisms. It is easy to show that

\begin{lem}
  Let $A=(A_0,A_1,A_2):\frkg\to \frkg'$ be
a homomorphism. If $A_0, A_1$ are invertible, then $A$ is an isomorphism, and $A^{-1}$ is given by
 $$A^{-1}=(A^{-1}_0, A^{-1}_1, -A^{-1}_1 \circ A_2\circ (A^{-1}_0\times A^{-1}_0)).$$
\end{lem}

To define $\Aut(\g)$, first denote by $\Aut^0(\g)$ the set of Lie $2$-algebra automorphisms\footnote{Strictly speaking, we should use the terminology ``a weak automorphism'' since it preserves the structure only up to homotopy. To be succinct, we just say ``an automorphism''.} of $\g$. It is evident that $(\Aut^0(\g),\diamond)$ is a Lie group.
Next, define a multiplication on $\End^{-1}(\g)$ by
\begin{eqnarray}\label{eq:multi 1}
\tau\star\tau'=\tau+\tau'+\tau\circ d\circ \tau',\quad \forall ~\tau,\tau' \in \End^{-1}(\g).
\end{eqnarray}
It is obvious that $\star$ satisfies the associative law. Thus, $(\End^{-1}(\g),\star)$ is a monoid, in which the zero map is the identity element.
$\Aut^{-1}(\g)$ is
defined to be the group of units of $\End^{-1}(\g)$, which is a Lie group.

\begin{lem}\label{lem:invertible}
For all $\tau\in \End^{-1}(\g)$, we have $$\tau \in \Aut^{-1}(\g)\Longleftrightarrow I+d\circ\tau \in GL(\g_0)\Longleftrightarrow I+\tau\circ d\in GL(\g_{-1}).$$
\end{lem}
\pf First we prove the equivalence of the first two formulas. Let $\tau \in \Aut^{-1}(\g)$, then there exists an element $\tau^{-1}\in \Aut^{-1}(\g)$ such that $\tau\star \tau^{-1}=\tau^{-1}\star \tau=0.$ We have \begin{equation}\label{eq:inverse}
  (I+d\circ\tau)^{-1}=I+d\circ\tau^{-1},
\end{equation} which follows from
\begin{eqnarray*}
(I+d\circ\tau)\circ (I+d\circ\tau^{-1})&=&I+d\circ(\tau+\tau^{-1}+\tau\circ  d\circ\tau^{-1})=I,\\
(I+d\circ\tau^{-1})\circ (I+d\circ\tau)&=&I+d\circ(\tau+\tau^{-1}+\tau^{-1}\circ d\circ\tau)=I.
\end{eqnarray*}
For the inverse direction, if the inverse of $I+d\circ\tau$ exists, then we claim that $\tau^{-1}=-\tau\circ (I+d\circ \tau)^{-1}.$
In fact, we have
\begin{eqnarray*}
\tau\star (-\tau\circ (I+d\circ \tau)^{-1})&=&\tau-\tau\circ(I+d\circ \tau)^{-1}-\tau\circ d\circ \tau\circ(I+d\circ \tau)^{-1}\\
&=&\tau-\tau\circ(I+d\circ \tau)\circ(I+d\circ \tau)^{-1}\\&=&0.
\end{eqnarray*}
Similarly, we have $(-\tau\circ (I+d\circ \tau)^{-1})\star \tau=0.$

The equivalence of the first and third formulas can be  proved similarly. The proof is finished.\qed\vspace{3mm}

Now, we are in the position to construct a crossed module from Lie groups $(\Aut^0(\g),\diamond)$ and $(\Aut^{-1}(\g),\star)$. The following observation plays an important role in the construction of the homomorphism from $\Aut^{-1}(\g)$ to $\Aut^0(\g)$.

\begin{lem}
Let $A:\g\longrightarrow\g$ be a Lie $2$-algebra homomorphism. Then for any $\tau\in\End^{-1}(\g)$,
$(A_0+d\circ \tau,A_1+\tau\circ d, A_2+l^A_\tau)$ is a Lie $2$-algebra homomorphism, where
$$l^A_\tau(x,y)=\tau[x,y]-[A_0x,\tau y]-[\tau x,A_0y]-[\tau x,d \tau y].$$
In particular, $(I+d\circ \tau,I+\tau\circ d, l^{\ie}_\tau)$  is a Lie $2$-algebra homomorphism, where $\ie=(I,I,0)$ is the identity homomorphism.
\end{lem}
\pf Firstly, due to $A_0\circ d= d\circ A_1,$ it is obvious that
\begin{equation}\label{eq:t1}
  (A_0+d\circ \tau)\circ d= d\circ (A_1+\tau\circ d).
\end{equation}
Then, by the fact that $A$ is a Lie $2$-algebra homomorphism, we have
\begin{eqnarray}
\nonumber&&(A_0+d\circ \tau)[x,y]-[(A_0+d\circ \tau)x,(A_0+d\circ \tau)y]\\ \nonumber&=&dA_2(x,y)+d\tau[x,y]-d[A_0x,\tau y]-d[\tau x,A_0y]-d[\tau x,d\tau y]\\\label{eq:t2} &=&d(A_2+l^A_\tau)(x,y).
\end{eqnarray}
Similarly, we can deduce that
\begin{equation}\label{eq:t3}(A_1+\tau\circ d)[x,a]-[(A_0+d\circ \tau)x,(A_1+\tau\circ d)a]=(A_2+l^A_\tau)(x,da).\end{equation}
Finally, also by the fact that $A$ is a Lie $2$-algebra homomorphism, we have
\begin{eqnarray}\label{cohe}
&&\nonumber [(A_0+d\circ \tau)x,(A_2+l^A_\tau)(y,z)]-(A_2+l^A_\tau)([x,y],z)+c.p.
\\ \nonumber&=&[A_0x,A_2(y,z)]-A_2([x,y],z)+c.p.\\\nonumber
&&+[d\tau x,A_2(y,z)]+[d\tau x+A_0x,\tau[y,z]-[A_0y,\tau z]-[\tau y,A_0z]-[\tau y,d\tau z]]\\\nonumber &&
-\tau[[x,y],z]+[[A_0x,A_0y]+dA_2(x,y),\tau z]+[\tau[x,y],A_0z+d\tau z]+c.p.\\ \nonumber
&=&A_1l_3(x,y,z)-l_{3}(A_0x,A_0y,A_0z)
+\tau dl_3(x,y,z)\\ \nonumber
&&-[d\tau x+A_0x,[A_0y,\tau z]+[\tau y,A_0z]+[\tau y,d\tau z]]+[[A_0x,A_0y],\tau z]+c.p.\\ \label{eq:t4}
&=&(A_1+\tau\circ d)l_3(x,y,z)-l_{3}\big((A_0+d\circ \tau)x,(A_0+d\circ \tau)y,(A_0+d\circ \tau)z\big).
\end{eqnarray}
By \eqref{eq:t1}-\eqref{eq:t4}, we deduce that $(A_0+d\circ \tau,A_1+\tau\circ d, A_2+l^A_\tau)$ is also a Lie $2$-algebra homomorphism. This ends the proof.\qed\vspace{3mm}


By Lemma \ref{lem:invertible}, if $\tau\in\Aut^{-1}(\g)$, then $(I+d\circ \tau,I+\tau\circ d, l^{\ie}_\tau)\in \Aut^0(\g)$. Thus, we have a smooth map $\partial:\Aut^{-1}(\g)\longrightarrow \Aut^0(\g)$, which is given by
\begin{equation}\label{eq:partial}
  \partial(\tau)=(I+d\circ \tau,I+\tau\circ d, l^{\ie}_\tau),\quad \forall \tau\in\Aut^{-1}(\g).
\end{equation}
Furthermore, $\Aut^0(\g)$ acts on  $\Aut^{-1}(\g)$ naturally:
\begin{equation}\label{eq:action}
 A\triangleright \tau=A_1\circ\tau\circ A^{-1}_0\quad \forall A=(A_0,A_1,A_2)\in\Aut^0(\g),~ \tau\in\Aut^{-1}(\g).
\end{equation}

\begin{thm}\label{thm:2grp}
With the notations above,  $(\Aut^{-1}(\g),\Aut^0(\g),\partial,\triangleright)$ is a crossed module of Lie groups, which is called the {\bf automorphism $2$-group} of $\g$,  and denoted by $\Aut(\g)$.
\end{thm}
\pf
For any $\tau, \tau'\in\Aut^{-1}(\g)$, we have
\begin{eqnarray*}
(\partial\tau)\diamond (\partial \tau')&=&(I+d\circ\tau,I+\tau\circ d,l^{\ie}_\tau)\diamond(I+d\circ\tau',I+\tau'\circ d,l^{\ie}_{\tau'})\\&=&(I+d\circ(\tau\star \tau'),I+(\tau\star \tau')\circ d,F_2),
\end{eqnarray*}
where $F_2$ is given by
\begin{eqnarray*}
F_2(x,y)&=&l^{\ie}_\tau((I+d\circ\tau')x,(I+d\circ\tau')y)+(I+\tau\circ d)l^{\ie}_{\tau'}(x,y)\\ &=&
\tau[x+d\tau'x,y+d\tau'y]-[x+d\tau'x,\tau(y+d\tau'y)]-[\tau(x+d\tau'x),y+d\tau'y]\\ &&-[\tau(x+d\tau'x)
,d\tau(y+d\tau'y)]+(I+\tau\circ d)(\tau'[x,y]-[\tau'x,y]-[x,\tau'y]-[\tau'x,d\tau'y])\\ &=&
\tau[x,y]+\tau'[x,y]+\tau d \tau'[x,y]-[\tau x+\tau'x+\tau d\tau'x,y]\\ &&-[x,\tau y+\tau'y+\tau d\tau'y]-[\tau x+\tau'x+\tau d\tau'x,d\tau y+d\tau'y+d\tau d\tau'y]\\&=&l^{\ie}_{\tau\star \tau'}(x,y).
\end{eqnarray*}
Thus, we have
 \begin{equation}\label{eq:pm}
  \partial(\tau\star \tau')=(\partial\tau)\diamond (\partial \tau'),
\end{equation}which implies that $\partial$ is a Lie group homomorphism.

For any $A\in \Aut^0(\g)$, we  have that $A^{-1}=(A^{-1}_0, A_1^{-1},-A_1^{-1}\circ A_2\circ (A^{-1}_0\times A^{-1}_0))$, and
\begin{eqnarray*}
A\diamond (\partial \tau)\diamond A^{-1}&=&(A_0,A_1,A_2)\diamond (I+d \circ \tau,I+\tau \circ d,l^{\ie}_\tau)\diamond (A^{-1}_0, A_1^{-1},-A_1^{-1}\circ A_2\circ (A^{-1}_0\times A^{-1}_0))\\ &=&\big(I+d\circ(A\triangleright\tau),I+(A\triangleright\tau)\circ d, H_2),
\end{eqnarray*}
where
\begin{eqnarray*}
&&H_2(x,y)\\ &=& \big(A\diamond (\partial \tau)\big)_2(A^{-1}_0x,A^{-1}_0y)-A_1\circ(I+\tau\circ d)\circ A_1^{-1}\circ A_2\circ (A^{-1}_0\times A^{-1}_0)(x,y)\\ &=&A_2\big((I+d\circ\tau)(A^{-1}_0x),(I+d\circ\tau)(A^{-1}_0y)\big)-A_2(A^{-1}_0x,A^{-1}_0y)-A_1\tau dA_1^{-1}A_2(A^{-1}_0x,A^{-1}_0y)\\ &&+A_1\big(\tau[A^{-1}_0x,A^{-1}_0y]-[\tau A^{-1}_0x, A^{-1}_0y]-[ A^{-1}_0x, \tau A^{-1}_0y]-[\tau A^{-1}_0x,d \tau A^{-1}_0y]\big)\\ &=&
A_1\tau A_0^{-1}[x,y]-[A_1\tau A_0^{-1}x,y]
-[x,A_1\tau A_0^{-1}y]-[A_1\tau A_0^{-1}x,dA_1\tau A_0^{-1}y]\\ &=&l^{\ie}_{A_1\circ \tau \circ A_0^{-1}}(x,y)=l_{A\triangleright\tau}^{\ie}(x,y).
\end{eqnarray*}
Thus, we have \begin{equation}\label{eq:c1}
  A\diamond (\partial \tau)\diamond A^{-1}=\partial(A\triangleright\tau).
\end{equation}

By straightforward computations, we have
\begin{eqnarray*}
\tau \star \tau'\star \tau^{-1}&=&\tau \star \tau' +\tau^{-1}+(\tau \star \tau')\circ d\circ \tau^{-1}\\ &=&
\tau +\tau'+\tau\circ d\circ \tau'+\tau^{-1}+\tau\circ d\circ\tau^{-1}+\tau'\circ d\circ\tau^{-1}+\tau\circ d\circ\tau'\circ d\circ\tau^{-1}\\ &=&
\tau'+\tau\circ d\circ\tau'+\tau'\circ d\circ\tau^{-1}+\tau\circ d\circ\tau'\circ d\circ\tau^{-1}.
\end{eqnarray*}
On the other hand,  by \eqref{eq:inverse}, we have
\begin{eqnarray*}
\partial(\tau)\triangleright \tau'&=&(I+\tau\circ d)\circ\tau'\circ(I+d\circ\tau)^{-1}=
(I+\tau \circ d)\circ\tau'\circ(I+d\circ \tau^{-1})\\ &=& \tau'+\tau\circ d\circ\tau'+\tau'\circ d\circ\tau^{-1}+\tau\circ d\circ\tau'\circ d\circ\tau^{-1}.
\end{eqnarray*}
Thus, we have
\begin{equation}\label{eq:c2}
  \partial(\tau)\triangleright \tau'=\tau \star \tau'\star \tau^{-1}.
\end{equation}
By \eqref{eq:c1} and \eqref{eq:c2}, we deduce that $(\Aut^{-1}(\g),\Aut^0(\g),\partial,\triangleright)$ is a crossed module of Lie groups. The proof is finished.\qed

\begin{ex}{\rm ({\bf the string Lie 2-algebra $\str(\frkk)$})
For the string Lie 2-algebra $\str(\frkk)$, we explore its automorphism $2$-group $
\Aut(\str(\frkk))$. For any $(A_0,t)\in\End(\frkk)\oplus\mathbb R$
and $A_2\in \wedge^2\frkk^*$, we have
 \begin{equation*}
A\in\Aut^0(\str(\frkk)) \Longleftrightarrow \left\{\begin{array}{rcll}A_0 &\in& \Aut(\frkk ),\\
\mathfrak{D} A_2(x,y,z)&=&tl_3(x,y,z)-l_3(A_0x,A_0y,A_0z).\end{array}\right.
\end{equation*}
Since $\frkk$ is semisimple, any automorphism is orthogonal. Then the second condition reduces to $$\mathfrak{D} A_2(x,y,z)=tl_3(x,y,z)-l_3(x,y,z).$$
Also by the fact that the Cartan 3-form is not exact, we deduce that $t=1$ and $\frkD (A_2)=0$. Thus, any  automorphism   $A_0\in\Aut(\frkk)$ and a 2-cocycle $\omega\in\wedge^2\frkk^*$ give rise to an automorphism $(A_0,1,\omega)$ of degree 0 of the string Lie 2-algebra. Parallel to the discussion in Proposition \ref{pro:string}, we obtain that $\Aut^0(\str(\frkk)) $ is isomorphic to the semidirect product $\huaK\ltimes\frkk^*$, where $\huaK$ is the connected and simply connected Lie group that integrating $\frkk.$

Obviously, the set of automorphisms of degree $-1$ is $\frkk^*$, which is a Lie group with the abelian group structure. The map $\partial:
\Aut^{-1}(\str(\frkk))\rightarrow \Aut^0(\str(\frkk))$ is given by
$\partial(\tau)=(I,I,-\mathfrak{D}(\tau))$.
  }
\end{ex}

\begin{ex}\rm{
For the skeletal Lie $2$-algebra $\g=\frkk\oplus V$ given in Example \ref{ex of der}, we explore its automorphism $2$-group $
\Aut(\g)$. Indeed, for any $(A_0,A_1)\in\End(\frkk)\oplus\End(V)$
and $A_2\in C^2(\frkk,V)$, it is easy to check that
\begin{equation*}
A\in\Aut^0(\g) \Longleftrightarrow \left\{\begin{array}{rcll}A_0\oplus A_1&\in& \Aut(\frkk\ltimes_\phi V),\\
\mathfrak{D}^A A_2&=&[A,l_3],\end{array}\right.
\end{equation*}
where $\mathfrak{D}^A: C^k(\frkk,V)\rightarrow C^{k+1}(\frkk,V)$ is the
Lie algebra coboundary operator with respect to the representation $\phi^A_x (v)=\phi_{A_0(x)}(v)$,  $\Aut(\frkk\ltimes_\phi V)$ is the Lie
group of automorphisms of $\frkk\ltimes_\phi
V$ and the bracket $[\cdot,\cdot]$ is given by
$$[A,l_3](x,y,z)=A_1l_3(x,y,z)-l_3(A_0x,A_0y,A_0z),\ \ \ \ \ \ \ \forall x,y,z\in \frkk.$$
Moreover, we
have $\Aut^{-1}(\g)=C^1(\frkk,V)$, which is a Lie group with the abelian group structure. The map $\partial:
\Aut^{-1}(\g)\rightarrow \Aut^0(\g)$ is given by
$\partial(\tau)=(I,I,-\mathfrak{D}(\tau))$.}
\end{ex}

\section{The integration of the Lie 2-algebra  $\Der(\g)$}
In this section, we show that for a Lie 2-algebra $\g$, the automorphism 2-group $\Aut(\g)$ is an integration of the derivation Lie 2-algebra $\Der(\g)$, i.e.  the differentiation of  $\Aut(\g)$ is the crossed module of Lie algebras corresponding to $\Der(\g)$. In addition, we give the expressions of exponential maps from $\Der^i(\g)$ to $\Aut^i(\g)$ for $i=0,-1$, which commute with the differentials $\partial$ and $\bar{d}$.

\begin{thm}\label{thm:main}
The differentiation of the crossed module of Lie groups $(\Aut^{-1}(\g),\Aut^0(\g),\partial,\triangleright)$ is the crossed module of Lie algebras $(\Der^{-1}(\g),\Der^0(\g),\bar{d},\phi)$.
\end{thm}
We split the proof into three steps.

$\bullet$ Step $1$: We shall check the Lie algebra of $\Aut^0(\g)$ is $\Der^0(\g)$. Define a map
$$\exp: \mathbb R\times \big(\End(\g_0)\oplus \End(\g_{-1})\oplus\Hom(\wedge^2\g_0,\g_{-1})\big)\longrightarrow\End(\g_0)\oplus \End(\g_{-1})\oplus\Hom(\wedge^2\g_0,\g_{-1})$$ by
\begin{eqnarray}\label{exponential1}
\exp(t,(X_0,X_1,l_X))=(e^{tX_0},e^{tX_1},e^{tl_X}),
\end{eqnarray}
 where $e^{tX_i}=\sum_{n=0}^{\infty}\frac{t^{n}X_i^n}{n!},i=0,1,$ is the usual exponential map, and $e^{tl_X}$ is given by
\begin{equation}\label{eq:elx}
e^{tl_X}(x,y)=\sum_{n=1}^{\infty}\frac{t^{n}}{n!}\sum_{i+j+k=n-1, i,j,k\geq 0}C_{i+j}^iX_1^kl_X(X_0^ix,X_0^jy),\ \ \ \ \forall x,y \in \g_0.
\end{equation}
Later, we will see that $\exp(1,(X,l_X)): \Der^0(\g)\longrightarrow \Aut^0(\g)$ is the exponential map.
\begin{lem}\label{lem:tangent space}
The tangent space of $\Aut^0(\g)$ at the identity element $\ie=(I_{\g_0},I_{\g_{-1}},0)$ is $\Der^0(\g)$, i.e. $$T_\ie \Aut^0(\g)=\Der^0(\g).$$
\end{lem}
\pf
Firstly, we prove $T_\ie \Aut^0(\g)\subset\Der^0(\g)$. Let $c(t)=(c_0(t),c_1(t),c_2(t))$ be any curve in $\Aut^0(\g)$ such that $c(0)=(I_{\g_0},I_{\g_{-1}},0)$ and $c'(0)=(X_0,X_1,l_X)$. By the fact that $c(t)$ is a curve in $\Aut^0(\g)$, we have
\begin{eqnarray*}
 && \frac{d}{dt}|_{t=0}\Big([c_0(t)x,c_2(t)(y,z)]-c_2(t)([x,y],z)+c.p.+
l_{3}(c_0(t)x,c_0(t)y,c_0(t)z)-c_1(t)l_{3}(x,y,z)\Big)\\
&=&[x,l_X(y,z)]-l_X([x,y],z)+l_3(X_0x,y,z)+c.p.-X_1l_3(x,y,z),
\end{eqnarray*}
which implies that one can obtain Condition (c) in Definition \ref{defi:derivation} from Condition (iii) in Definition \ref{defi:Lie 2 homo} by differentiation.
Similarly, we can obtain (a)-(b) from (i)-(ii) by differentiation. Thus,
 $c'(0)\in\Der^0(\g)$.
This proves that $T_\ie \Aut^0(\g)\subset\Der^0(\g)$.

On the other hand, for all $(X,l_X)\in \Der^0(\g)$,  we first claim that
$$c(t):=\exp(t,(X,l_X))=(e^{tX_0},e^{tX_1},e^{tl_X})\in \Aut^0(\g),\quad \forall t\in \mathbb R.$$
It is obvious that $e^{tX_0}$ and $e^{tX_1}$ are invertible and commute with $d$ since $X_0\circ d=d \circ X_1$. Consider the following two curves in $\g_0$: $$y_1(t)=e^{tX_0}[x,y],\ \ \ \ \ \ \ \ \ y_2(t)=[e^{tX_0}x,e^{tX_0}y]+de^{tl_X}(x,y).$$ Obviously, we have
\begin{equation}\label{eq:e1}
y_1(0)=y_2(0)=[x,y],
\end{equation}
and $y_1(t) $ satisfies the following ordinary differential equation:
\begin{equation}\label{eq:e2}
y'_1(t)=X_0e^{tX_0}[x,y]=X_0y_1(t).
\end{equation}
For any $t$, we have
\begin{eqnarray*}\label{diff}
\frac{d}{dt}\big(e^{tl_X}(x,y))\big)\nonumber&=&
\sum_{n=1}^{\infty}\frac{t^{n-1}}{(n-1)!}\sum_{i+j+k=n-1}C_{i+j}^iX_1^kl_X(X_0^ix,X_0^jy)\\ \nonumber&=&
\sum_{n=1}^{\infty}\frac{t^{n-1}}{(n-1)!}\sum_{i+j=n-1}\frac{(n-1)!}{i!j!}l_X(X_0^ix,X_0^jy)\\ \nonumber&&+
X_1\sum_{n=2}^{\infty}\frac{t^{n-1}}{(n-1)!}\sum_{i+j+k-1=n-2}C_{i+j}^iX_1^{k-1}l_X(X_0^ix,X_0^jy)
\\&=&l_X(e^{tX_0}x,e^{tX_0}y)+X_1e^{tl_X}(x,y).
\end{eqnarray*}
Hence, we get
\begin{eqnarray}
\nonumber y'_2(t)&=&[X_0e^{tX_0}x,e^{tX_0}y]+[e^{tX_0}x,X_0e^{tX_0}y]+dl_X(e^{tX_0}x,e^{tX_0}y)+X_0de^{tl_X}(x,y)\\ &=&\nonumber
X_0[e^{tX_0}x,e^{tX_0}y]+X_0de^{tl_X}(x,y)\\\label{eq:e3}&=&X_0y_2(t),
\end{eqnarray}
where we have used the fact
$$X_0[x,y]=[X_0x,y]+[x,X_0y]+dl_X(x,y),\ \ \ \ \ \ \ \forall x,y\in \g_0.$$
By \eqref{eq:e1}-\eqref{eq:e3} and the uniqueness theorem for linear systems
of ordinary differential equations, we obtain that $y_1(t)=y_2(t)$, i.e.
\begin{eqnarray}\label{eq:ee1}
e^{tX_0}[x,y]=[e^{tX_0}x,e^{tX_0}y]+
de^{tl_X}(x,y).
\end{eqnarray}
Similarly, we can get
\begin{eqnarray}\label{eq:ee2}
e^{tX_1}[x,a]=[e^{tX_0}x,e^{tX_1}a]+
e^{tl_X}(x,da).
\end{eqnarray}
At last, consider the following two curves in $\g_{-1}$:
\begin{eqnarray*}
 z_1(t)&=&-e^{tX_1}l_3(x,y,z),\\
  z_2(t)&=&e^{tl_X}([x,y],z)+[e^{tl_X}(x,y),e^{tX_0}z]+c.p.-l_3(e^{tX_0}x,e^{tX_0}y,e^{tX_0}z).
\end{eqnarray*}
Obviously, we have
$z_1(0)=z_2(0)=-l_3(x,y,z),$
and $z_1(t) $ satisfies the  ordinary differential equation:
$
z'_1(t)=X_1z_1(t).
$
By the fact $(X,l_X)\in \Der^0(\g)$, we can deduce that $z'_2(t)=X_1z_2(t)$. Thus, we have $z_2(t)=z_1(t)$, i.e.
\begin{equation}\label{eq:ee3}
e^{tl_X}([x,y],z)+[e^{tl_X}(x,y),e^{tX_0}z]+c.p.-l_3(e^{tX_0}x,e^{tX_0}y,e^{tX_0}z)=-e^{tX_1}l_3(x,y,z).
\end{equation}
By \eqref{eq:ee1}-\eqref{eq:ee3}, we deduce that $(e^{tX_0},e^{tX_1},e^{tl_X})\in \Aut^0(\g)$ for all $t\in \mathbb R$.

Therefore, for all $(X,l_X)\in \Der^0(\g)$, we get a curve $c(t):=(e^{tX_0},e^{tX_1},e^{tl_X})$ in $\Aut^0(\g)$ satisfying
$c(0)=\ie=(I_{\g_0},I_{\g_{-1}},0)$. It is obvious that $c'(0)=(X,l_X)$, which implies that $\Der^0(\g)\subset T_\ie \Aut^0(\g)$. This ends the proof. \qed\vspace{3mm}

For all $(X,l_X)\in \Der^0(\g)=T_\ie \Aut^0(\g)$, we get a smooth map
\begin{eqnarray*}
\sigma_{(X,l_X)}: \mathbb R\longrightarrow\Aut^0(\g)
\end{eqnarray*}
defined by $$\sigma_{(X,l_X)}(t)=\exp(t,(X,l_X))=(e^{tX_0},e^{tX_1},e^{tl_X}).$$
Define a map $e:\Der^0(\g)\longrightarrow \Aut^0(\g)$ by
\begin{equation}\label{eq:defi e}
  e^{(X,l_X)}=\sigma_{(X,l_X)}(1).
\end{equation}

\begin{lem}\label{lem:exp}
With the notations above, for all $(X,l_X)\in \Der^0(\g)$, the map $\sigma_{(X,l_X)}(t)$ is a one-parameter subgroup of the Lie group $\Aut^0(\g)$ determined by $(X,l_X)$. Consequently, the map $e:\Der^0(\g)\longrightarrow \Aut^0(\g)$ given by \eqref{eq:defi e} is the exponential map.
\end{lem}
\pf For simplicity, we omit the subscript $(X,l_X)$ when there is no confusion. It is obvious that
$e^{(t+s)X_i}=e^{tX_i}\circ e^{sX_i}.$
Furthermore, by the equality $ C_{i+j}^{i}=\sum_{\alpha=0}^{m}C_m^\alpha C_{i+j-m}^{i-\alpha},$
for any $i+j\geq m$, we have
\begin{eqnarray*}
e^{(t+s)l_X}(x,y) &=&\sum_{n=1}^{\infty}\frac{(t+s)^{n}}{n!}\sum_{i+j+k=n-1}C_{i+j}^{i}X_1^kl_X(X_0^ix,X_0^jy)\\ &=&
\sum_{n=1}^{\infty}\sum_{m=0}^{n}\frac{s^mt^{n-m}}{m!(n-m)!}\sum_{i+j+k=n-1}C_{i+j}^{i}X_1^kl_X(X_0^ix,X_0^jy)\\
&=&\sum_{p,m=0,p+m\geq 1}^{\infty}\frac{t^{p}}{p!}\frac{s^m}{m!}
\sum_{i+j+k=p+m-1}C_{i+j}^{i}X_1^k l_X(X_0^ix,X_0^jy)\\&=&\sum_{p=0}^{\infty}\frac{t^{p}}{p!}\sum_{m=1}^{\infty}\frac{s^m}{m!}
\sum_{i+j+k-p=m-1}C_{i+j}^{i}X_1^p X_1^{k-p}l_X(X_0^ix,X_0^jy)\\ &&+\sum_{p=1}^{\infty}\frac{t^{p}}{p!}\sum_{m=0}^{\infty}
\sum_{u+v+k=p-1}\sum_{\alpha=0}^{m}\frac{s^{m}C_{u+v}^{u}}
{\alpha!(m-\alpha)!}X_1^{k}l_X(X_0^{u} X_0^{\alpha}x,X_0^{v}X_0^{m-\alpha}y)
\\ &=&e^{tX_1}e^{sl_X}(x,y)+e^{tl_X}(e^{sX_0}x,e^{sX_0}y)\\ &=&\big((e^{tX_0},e^{tX_1},e^{tl_X})\diamond(e^{sX_0},e^{sX_1},e^{sl_X})\big)_2(x,y).
\end{eqnarray*}
Therefore, we have
$$(e^{(t+s)X_0},e^{(t+s)X_1},e^{(t+s)l_X})=(e^{tX_0},e^{tX_1},e^{tl_X})\diamond(e^{sX_0},e^{sX_1},e^{sl_X}),
$$
which implies that $\sigma(t+s)=\sigma(t)\diamond\sigma(s)$, and $\sigma(t)$ is a one-parameter subgroup.

It is obvious that $$\sigma'(0)=\frac{d}{dt}|_{{t=0}}(e^{tX_0},e^{tX_1},e^{tl_X})=(X,l_X),$$ which implies that $\sigma(t)$ is a one-parameter subgroup determined by $(X,l_X)$. Thus, the map $e:\Der^0(\g)\longrightarrow \Aut^0(\g)$ given by \eqref{eq:defi e} is the exponential map.
\qed\vspace{3mm}

With regard to the relation between the Lie group structure on $K$ and the Lie bracket of its Lie algebra $\frkk$, there is a useful formula:
\begin{eqnarray}\label{Lie alg of Lie gp}
e^{sX} e^{tY} e^{-sX} e^{-tY}\doteq e^{st[X,Y]_{\frkk}},\ \ \ \ \ \ \ \ \forall X,Y\in \frkk,
\end{eqnarray}
where $e:\frkk\longrightarrow K$ is the exponential map.

\begin{lem}\label{pro:int degree0}
The Lie algebra of the Lie group $(\Aut^0(\g),\diamond)$ is  $(\Der^0(\g),\{\cdot,\cdot\})$, where  the bracket $\{\cdot,\cdot\}$ is given by \eqref{bra of der}.
\end{lem}
 \pf By Lemma \ref{lem:tangent space}, $\Der^0(\g)$ is the tangent space of $\Aut^0(\g)$ at the identity element. To prove that the Lie algebra of the Lie group $(\Aut^0(\g),\diamond)$ is  $(\Der^0(\g),\{\cdot,\cdot\})$, we only need to show that the induced Lie bracket $[\cdot,\cdot]_{ind}$ on $\Der^0(\g)$ is exactly $\{\cdot,\cdot\}$. By Lemma \ref{lem:exp}, we have the exponential map $e:\Der^0(\g)\longrightarrow \Aut^0(\g)$ given by \eqref{eq:defi e}. Thus, we have
\begin{eqnarray*}
&&[(X,l_X),(Y,l_Y)]_{ind}\\ &=&\frac{d}{dt}\frac{d}{ds}|{_{t,s=0}}(e^{sX_0},e^{sX_1},e^{sl_X}) \diamond (e^{tY_0},e^{tY_1},e^{tl_Y})\diamond (e^{-sX_0},e^{-sX_1},e^{-sl_X})\diamond (e^{-tY_0},e^{-tY_1},e^{-tl_Y})\\ &=&
\frac{d}{dt}\frac{d}{ds}|{_{t,s=0}}(e^{sX_0}e^{tY_0}e^{-sX_0}e^{-tY_0},e^{sX_1}e^{tY_1}e^{-sX_1}e^{-tY_1},F_2)\\ &=&
\big([X_0,Y_0]_C,[X_1,Y_1]_C,\frac{d}{dt}\frac{d}{ds}{|_{t,s=0}}F_2\big),
\end{eqnarray*}
where  $F_2$ is given by
\begin{eqnarray*}
F_2(x,y)&=&
e^{sl_X}(e^{tY_0}e^{-sX_0}e^{-tY_0}x,e^{tY_0}e^{-sX_0}e^{-tY_0}y)+e^{sX_1}e^{tl_Y}(e^{-sX_0}e^{-tY_0}x,e^{-sX_0}e^{-tY_0}y)
\\ &&+e^{sX_1}e^{tY_1}e^{-sl_X}(e^{-tY_0}x,e^{-tY_0}y)+e^{sX_1}e^{tY_1}e^{-sX_1}e^{-tl_Y}(x,y).
\end{eqnarray*}
By  (\ref{exponential1}), we get
\begin{eqnarray*}
\frac{d}{dt}\frac{d}{ds}|{_{t,s=0}}F_2&=&X_1l_Y(x,y)-l_Y(X_0x,y)-l_Y(x,X_0y)-Y_1l_X(x,y)+l_X(Y_0x,y)+l_X(x,Y_0y)
\\ &=&L_Xl_Y(x,y)-L_Yl_X(x,y).
\end{eqnarray*}
So the induced Lie bracket is
$$[(X,l_X),(Y,l_Y)]_{ind}=([X,Y]_C,L_X l_Y-L_Y l_X),$$
which is exactly the Lie bracket (\ref{bra of der}) on $\Der^0(\g)$. This ends the proof.
\qed\vspace{3mm}

$\bullet$ Step $2$: We verify that the Lie algebra of $\Aut^{-1}(\g)$ is $\Der^{-1}(\g)$. Moreover, the two exponential maps from $\Der^i(\g)$ to $\Aut^i(\g)$ for $i=0,-1$ commute with the differentials $\partial$ and $\bar{d}$.
\begin{lem}\label{pro:int degree1}
The Lie algebra of $(\Aut^{-1}(\g),\star)$ is $(\Der^{-1}(\g),\{\cdot,\cdot\})$, where the multiplication $\star$ and the bracket $\{\cdot,\cdot\}$ are given by \eqref{eq:multi 1} and \eqref{eq:bra of der0} respectively. Additionally, the exponential map $e:\Der^{-1}(\g)\longrightarrow\Aut^{-1}(\g)$ is given by
\begin{eqnarray}\label{eq:defi e2}
e^\Theta=\Theta+\frac{\Theta\circ d\circ \Theta}{2!}+\frac{\Theta\circ d\circ \Theta\circ d \circ \Theta}{3!}+\cdots.
\end{eqnarray}
\end{lem}
\pf For all $\Theta\in \Der^{-1}(\g)$, define a one-parameter map $\sigma_\Theta:\mathbb R\longrightarrow \End^{-1}(\g)$ by
$\sigma_\Theta(t)=e^{t\Theta}.$
Then, we have  $I+d\circ \sigma_\Theta(t)=e^{t(d\circ \Theta)}\in GL(\g_0).$ By Lemma \ref{lem:invertible},
 $\sigma_\Theta(t)\in \Aut^{-1}(\g)$. It is straightforward to verify that $T_0\Aut^{-1}(\g)=\Der^{-1}(\g)$. Furthermore, $\sigma_\Theta$ is a one-parameter subgroup of $\Aut^{-1}(\g)$ determined by $\Theta$, which follows from
\begin{eqnarray*}
\sigma_\Theta(t+s)&=&\sum_{n=1}^{\infty}\sum_{m=0}^n\frac{t^m s^{n-m}}{m!(n-m)!}(\Theta\circ d)^{n-1}\circ \Theta\\ &=&
\sum_{n=1}^{\infty}\frac{ s^{n}}{n!}(\Theta\circ d)^{n-1}\circ \Theta+\sum_{n=1}^{\infty}\frac{t^n }{n!}(\Theta\circ d)^{n-1}\circ \Theta\\ &&+\sum_{n=1}^{\infty}\sum_{m=1}^{n-1}\frac{t^m s^{n-m}}{m!(n-m)!}(\Theta\circ d)^{m-1}\circ \Theta\circ d\circ (\Theta\circ d)^{n-m-1}\circ \Theta\\ &=&\sigma_\Theta(t)+\sigma_\Theta(s)+\sigma_\Theta(t)\circ d\circ \sigma_\Theta(s)\\ &=&\sigma_\Theta(t)\star\sigma_\Theta(s).
\end{eqnarray*}
So we get that $e:\Der^{-1}(\g)\longrightarrow \Aut^{-1}(\g)$ defined by (\ref{eq:defi e2}) is the exponential map determined by $\Theta$. 
 By \eqref{eq:bra of der0}, \eqref{eq:multi 1} and \eqref{Lie alg of Lie gp}, the induced Lie bracket on $\Der^{-1}(\g)$ by the Lie group structure on $\Aut^{-1}(\g)$
is given by
\begin{eqnarray*}
&&[\Theta,\Theta']_{ind}\\ &=&\frac{d}{dt}\frac{d}{ds}|{_{t,s=0}}e^{s\Theta} \star e^{t\Theta'}\star  e^{-s\Theta}\star e^{-t\Theta'}\\ &=&\frac{d}{dt}\frac{d}{ds}|{_{t,s=0}}(s\Theta+\cdots)\star(t\Theta'+\cdots)
\star(-s\Theta+\cdots)\star(-t\Theta'+\cdots)\\ &=&
\frac{d}{dt}\frac{d}{ds}|{_{t,s=0}}(s\Theta+t\Theta'+st\Theta d\Theta'+\cdots)\star(-s\Theta+\cdots)\star(-t\Theta'+\cdots)\\ 
\\ &=&\frac{d}{dt}\frac{d}{ds}|{_{t,s=0}}(st\Theta\circ d\circ \Theta'-ts\Theta'\circ  d\circ \Theta+\cdots)
\\ &=&\Theta\circ d\circ \Theta'-\Theta'\circ  d\circ \Theta\\
&=&\{\Theta,\Theta'\}.
\end{eqnarray*}
Thus, the Lie algebra of $(\Aut^{-1}(\g),\star)$ is $(\Der^{-1}(\g),\{\cdot,\cdot\})$. This completes the proof.\qed\vspace{3mm}

The two exponential maps given by \eqref{eq:defi e} and \eqref{eq:defi e2} commute with the differentials $\partial $ and $\overline{d}$.

\begin{lem}\label{lem:commu e}
For any $\Theta\in \Der^{-1}(\g)$,  $\partial(e^\Theta)=e^{\bar{d}(\Theta)},$ i.e. the following diagram is commutative:
$$
\xymatrix{
\Der^{-1}(\g)\ar[r]^{e}\ar[d]^{\overline{d}}&\Aut^{-1}(\g)\ar[d]^{\partial}\\
\Der^0(\g)\ar[r]^{e}&\Aut^{-1}(\g).}
$$
\end{lem}
\pf First, by definition, we have
$$I+d\circ e^\Theta=I+\sum_{n=1}^\infty\frac{d\circ(\Theta\circ d)^{n-1}\circ \Theta}{n!}=e^{d\circ \Theta},$$
and $I+e^\Theta\circ d=e^{\Theta\circ d}$. Then, by straightforward computation, we get
\begin{eqnarray*}
&&e^{l_{\delta(\Theta)}}(x,y)\\ &=&\sum_{n=1}^{\infty}\frac{1}{n!}\sum_{i+j+k=n-1}C_{i+j}^i(\Theta d)^kl_{\delta(\Theta)}((d\Theta)^ix,(d\Theta)^jy)\\ &=&\sum_{n=1}^{\infty}\frac{1}{n!}\sum_{i+j+k=n-1}C_{i+j}^i(\Theta d)^k(\Theta[(d\Theta)^ix,(d\Theta)^jy]-[\Theta(d\Theta)^ix,(d\Theta)^jy]-[(d\Theta)^ix,\Theta(d\Theta)^jy])\\ &=&
e^\Theta[x,y]-[e^\Theta x,y]-[x,e^\Theta y]+\sum_{n=1}^{\infty}\frac{1}{n!}\{\sum_{i'+j+k'=n-1,k'\geq1}C_{i'+j+1}^{i'+1}(\Theta d)^{k'}[(\Theta d)^{i'}\Theta x,(d\Theta)^jy]\\ &&-\sum_{i+j+k=n-1,j+k\geq1}C_{i+j}^i(\Theta d)^k[(\Theta d)^i\Theta x,(d\Theta)^jy]\\ &&-\sum_{i'+j'+k=n-1,j'\geq 1}C_{i'+j'}^{i'+1}(\Theta d)^k[(\Theta d)^{i'}\Theta x,(d\Theta)^{j'}y]\}\\ &=&
e^\Theta[x,y]-[e^\Theta x,y]-[x,e^\Theta y]-\sum_{n=1}^{\infty}\frac{1}{n!}\{\sum_{i'+j=n,i',j\geq1}C_{n-1}^{i'-1}[(\Theta d)^{i'-1}\Theta x,d(\Theta d)^{j-1}\Theta y]\\ &&+\sum_{i''+j'=n,i'',j'\geq1}C_{n-1}^{i''}[(\Theta d)^{i''-1}\Theta x,d(\Theta d)^{j'-1}\Theta y]\}\\ &=&
e^\Theta[x,y]-[e^\Theta x,y]-[x,e^\Theta y]-[e^\Theta x,d e^\Theta y]\\ &=&l_{e^\Theta}^{\ie}(x,y),
\end{eqnarray*}
where we have used the equalities $C_{p+q+1}^{p+1}=C_{p+q}^{p}+C_{p+q}^{p+1}$ and $C_{n-1}^{p-1}+C_{n-1}^p=C_{n}^p  $ in the third and fourth steps respectively. Therefore, we have $$\partial(e^\Theta)=(I+d\circ e^\Theta,I+e^\Theta\circ d,l_{e^\Theta}^{\ie})=(e^{d \circ \Theta},e^{\Theta \circ d},e^{l_{\delta(\Theta)}})=e^{\bar{d}(\Theta)}.$$ This finishes the proof.\qed\vspace{3mm}

$\bullet$ Step $3$: The final step is  to prove that the differential of $\partial$ is $\bar{d}$ and the induced action of $\triangleright$ on  $\Der(\g)$ is $\phi$.
By Lemma \ref{lem:commu e}, we have  $\partial(e^{t\Theta})=e^{t\bar{d}(\Theta)}$, which implies $\partial_*=\bar{d}$.
\emptycomment{
By the fact that $e^{t\Theta}{|_{t=0}}=0$, we have
\begin{eqnarray*}
\partial_*(\Theta)&=&\frac{d}{dt}{|_{t=0}}\partial(e^{t\Theta})=\frac{d}{dt}{|_{t=0}}
\big(I+d\circ e^{t\Theta},I+e^{t\Theta}\circ d,l_{e^{t\Theta}}^I\big)\\ &=&
\big(d\circ\Theta,\Theta\circ d,\frac{d}{dt}{|_{t=0}}\big(e^{t\Theta}[x,y]-[e^{t\Theta}x,y]-[x,e^{t\Theta}y]-[e^{t\Theta}x,de^{t\Theta}y]\big)\big)\\
&=&\big(d\circ\Theta,\Theta\circ d,\Theta[x,y]-[\Theta x,y]-[x,\Theta y]\big)\\
&=&\bar{d}\Theta.
\end{eqnarray*}
}

Denote by $\widetilde{\phi}$ the induced action of $\Der^0(\g)$ on $\Der^{-1}(\g)$. By straightforward calculations, we have
\begin{eqnarray*}
\widetilde{\phi}_{(X,l_X)}\Theta&=&\frac{d}{dt}\frac{d}{ds}|_{{t,s=0}}(e^{tX_0},e^{tX_1},e^{tl_X})\triangleright e^{s\Theta}=\frac{d}{dt}\frac{d}{ds}|_{{t,s=0}}e^{tX_1}\circ e^{s\Theta} \circ e^{-tX_0}\\ &=&
X_1\circ \Theta-\Theta\circ X_0=\{X+l_X,\Theta\}\\&=&\phi_{(X,l_X)}\Theta.
\end{eqnarray*}
This ends the proof of Theorem \ref{thm:main}. \qed
\vspace{3mm}

Consider the Lie algebra $\Der^0(\g)\ltimes\Der^{-1}(\g)$ with brackets defined by (\ref{eq:bra01})-(\ref{eq:bra of der0}) and the Lie group $\Aut^0(\g)\ltimes\Aut^{-1}(\g)$ with group structure (\ref{m h}) defined by
$$(A,\tau)\cdot(A',\tau')=(A\diamond A',\tau\star (A\triangleright \tau')).$$

\begin{cor} \label{cor:e}
The Lie algebra of $\Aut^0(\g)\ltimes\Aut^{-1}(\g)$ is $\Der^0(\g)\ltimes\Der^{-1}(\g)$. Moreover, the exponential map $e:\Der^0(\g)\ltimes\Der^{-1}(\g)\longrightarrow \Aut^0(\g)\ltimes\Aut^{-1}(\g)$ is given by
$$
e^{((X,l_X),\Theta)}=(e^{(X,l_X)},e^\Theta),
$$
where $e^{(X,l_X)}, e^\Theta$ are described explicitly by \eqref{eq:defi e} and \eqref{eq:defi e2} respectively.
\end{cor}

\section{Strict automorphisms and Inner automorphisms}
\subsection{Strict automorphisms}

In this part, we first consider a sub-Lie $2$-algebra $\SDer(\g)$ of $\Der(\g)$ and a sub-crossed module $\SAut(\g)$ of $\Aut(\g)$.
Then, for a strict Lie $2$-algebra $\g$ seen as a crossed module with the corresponding crossed module of Lie groups $\mathbb{G}=(G_1,G_0,d^{\mathbb G},\triangleright)$, we derive the relation between $\Aut(\mathbb G)$ introduced in \cite{Norrie} and $\SAut(\g)$.

For a Lie $2$-algebra $\g$, consider a sub-complex $\SDer(\g)$ of $\Der(\g)$ defined by
\begin{eqnarray*}
\SDer^0(\g)&\triangleq&\{X\in \End(\g_0)\oplus\End(\g_{-1})| (X,0)\in \Der^0(\g)\};\\
\SDer^{-1}(\g)&\triangleq&\{\Theta\in \Der^{-1}(\g)=\Hom(\g_0,\g_{-1})|\Theta[x,y]=[x,\Theta y]+[\Theta x,y],\forall x,y\in \g_0\}.
\end{eqnarray*}
We call them {\bf strict derivations} of a Lie 2-algebra. Explicitly, $X=(X_0,X_1)\in \SDer^0(\g)$ if and only if for any $x,y,z\in \g_0,a\in \g_{-1}$,
\begin{equation*}
\left\{\begin{array}{rcll} X_0\circ d&=&d\circ X_1,\\
X_0[x,y]&=&[X_0x,y]+[x,X_0y],\\
X_1[x,a]&=&[X_0x,a]+[x,X_1a],\\
X_1l_{3}(x,y,z)&=&l_{3}(X_0x,y,z)+c.p.(x,y,z).\end{array}\right.
\end{equation*}
Then consider a sub-complex $\SAut(\g)$ of $\Aut(\g)$ given by
\begin{eqnarray*}
\SAut^0(\g)&\triangleq&\{A\in \End(\g_0)\oplus\End(\g_{-1})|(A_0,A_1,0)\in \Aut^0(\g)\};\\
\SAut^{-1}(\g)&\triangleq&\{\tau \in \Aut^{-1}(\g)=\Hom(\g_0,\g_{-1})|\tau[x,y]=[x,\tau y]+[\tau x,y]+[\tau x,d\tau y],\forall x,y\in \g_0\}.
\end{eqnarray*}
In detail, $A=(A_0,A_1)\in \SAut^0(\g)$ if and only if $A_0\circ d=d\circ A_1,$ and
$$A_0[x,y]=[A_0x,A_0y],\ \ \ A_1[x,a]=[A_0x,A_1a],\ \ \ A_1l_{3}(x,y,z)=l_{3}(A_0x,A_0y,A_0z).$$
\begin{lem}
$(\SAut^{-1}(\g),\SAut^0(\g),\partial,\triangleright)$ is a sub-crossed module of Lie groups of $\Aut(\g)$. Moreover,
the differentiation of $\SAut(\g)$ is the crossed module of Lie algebras $(\SDer^{-1}(\g),\SDer^0(\g),\bar{d},\phi)$, which is a sub-crossed module of $\Der(\g)$.
\end{lem}
\pf It is a straightforward verification.\qed\vspace{3mm}

Next, suppose that $\g$ is a {\bf strict} Lie $2$-algebra with the corresponding crossed module of Lie groups $\mathbb{G}=(G_1,G_0,d^{\mathbb G},\triangleright)$.
According to \cite{Norrie}, the automorphism $2$-group $\Aut(\mathbb{G})$ of $\mathbb{G}$ is defined as follows.
Define $$\Der(G_0,G_1)=\{\chi:G_0\longrightarrow G_1|\chi(\alpha\beta)=\chi(\alpha)(\alpha\triangleright\chi(\beta)), \forall \alpha,\beta\in G_0 \},$$ with a multiplication
\begin{eqnarray}\label{ast}
\chi\ast\chi'(\alpha)=\chi(d^\mathbb{G}\chi'(\alpha)\alpha)\chi'(\alpha).
\end{eqnarray}
Denote by $\Aut^{-1}(\mathbb{G})$ the group of units of $\Der(G_0,G_1)$.
Define
$$\Aut^0(\mathbb{G}):=\{(F_0,F_1)\in \Aut(G_0)\times\Aut(G_1)|F_0\circ d^\mathbb{G}=d^\mathbb{G}\circ F_1, F_1(\alpha\triangleright\xi)=F_0(\alpha)\triangleright F_1(\xi)\}.$$
Indeed, $\Aut^0(\mathbb{G})$ is the set of isomorphisms of $\mathbb{G}$. Then $\Aut(\mathbb{G})$ becomes a crossed module of Lie groups
\begin{equation*}\CD \Aut(\mathbb{G}): \Aut^{-1}(\mathbb{G}) @>\tilde{\partial}>>\Aut^0(\mathbb{G}),\endCD\end{equation*}where $\tilde{\partial}$ is defined by
\begin{eqnarray}\label{1-ope}
\tilde{\partial}(\chi)(\alpha)=d^\mathbb{G}(\chi(\alpha))\alpha,\ \ \ \tilde{\partial}(\chi)(\xi)=\chi(d^\mathbb{G}(\xi))\xi,
\end{eqnarray}
and the action of $\Aut^0(\mathbb{G})$ on $\Aut^{-1}(\mathbb{G})$ is defined by $$(F_0,F_1)\triangleright\chi=F_1\circ \chi \circ F_0^{-1}.$$

\begin{lem}\label{lem:t3}
For any $F\in \Aut^0(\mathbb{G})$, let $F_{*}=(F_{0*},F_{1*})$ be the tangent map at the identity. Then,  $F_{*}\in \SAut^0(\g)$, 
Furthermore, we have $$(F\circ F')_{*}=F_{*}\diamond F'_{*},\ \ \ \ \ \forall F,F'\in \Aut^0(\mathbb{G}).$$
\end{lem}

\begin{lem}\label{lem:t1}
For any $\chi\in \Aut^{-1}(\mathbb{G})$, let $\chi_{*}$ be the tangent map at the identity. Then, $\chi_{*}\in \SAut^{-1}(\g)$.
Moveover, we have
$$(\chi\ast \chi')_{*}=\chi_{*}\star \chi'_{*},\ \ \ \ \ \forall \chi,\chi'\in \Aut^{-1}(\mathbb{G}).$$
where $\ast$ and $\star$ are defined by \eqref{ast} and \eqref{eq:multi 1} respectively.
\end{lem}
\pf The condition $\chi(\alpha\beta)=\chi(\alpha)\alpha\triangleright\chi(\beta)$ implies that
$(I_{G_0},\chi): G_0\ltimes G_1\longrightarrow G_0\ltimes G_1$, which is given by
$$
(I_{G_0},\chi)(\alpha,\xi)=(\alpha,\chi(\alpha)),
$$ is a Lie group homomorphism.
By the fact that the Lie algebra of $G_0\ltimes G_1$ is $\g_0\ltimes\g_{-1}$, we obtain that $(I,\chi_{*}):\g_0\ltimes\g_{-1}\longrightarrow \g_0\ltimes\g_{-1}$ is a Lie algebra homomorphism, which  implies that $$\chi_{*}[x,y]=[x,\chi_{*}y]+[\chi_{*}x,y]+[d\chi_{*}x,\chi_{*}y],$$ i.e. $\chi_{*}\in \SAut^{-1}(\g)$. The other conclusion is obvious.
\qed\vspace{3mm}

Moveover, by simple calculations, we have
\begin{lem}\label{lem:t2}
For all $F\in \Aut^0(\mathbb{G})$ and $\chi\in \Aut^{-1}(\mathbb{G})$, we have
$$(\tilde{\partial}(\chi))_{*}=\partial(\chi_{*}),\ \ \ \ \ \ (F\triangleright\chi)_{*}=F_{*}\triangleright \chi_{*},$$
where $\bar{\partial}$ and $\partial$ are defined by \eqref{1-ope} and \eqref{eq:partial} respectively.
\end{lem}

In summary, by Lemma \ref{lem:t3}-\ref{lem:t2}, we have
\begin{thm}
The map from $\Aut(\mathbb G)$ to $\SAut(\g)$ given by
$$
(F,\chi)\longmapsto(F_{*},\chi_{*}),\quad \quad \forall F\in\Aut^0(\mathbb G), \chi\in\Aut^{-1}(\mathbb G),
$$
is a homomorphism between crossed modules of Lie groups.
\end{thm}

\emptycomment{Depending on the discussion above, we say that $\Aut(\mathbb{G})$ is the integration of $S\Aut(\g)$.

There exists the natural adjoint action $AD:\mathbb{G}\longrightarrow\Aut(\mathbb{G})$, which is a morphism of crossed modules defined by
\begin{equation*}
\left\{\begin{array}{rcl} AD_0(\alpha)(\beta)&=&\alpha\beta\alpha^{-1},
\\ AD_0(\alpha)(\xi)&=&\alpha\xi\alpha^{-1}=\alpha\triangleright\xi,\\
AD_1(\xi)(\alpha)&=&Pr_{G_1}(\xi\alpha\xi^{-1})=\xi\alpha\triangleright\xi^{-1}.\end{array}\right.
\end{equation*}
By taking differentiation, we get a map $ad: \g\longrightarrow S\Der(\g)$ defined as£»
\begin{equation*}
\left\{\begin{array}{rcl} ad_0(x)&=&[x,\cdot],
\\
ad_1(a)&=&[a,\cdot],\end{array}\right.
\end{equation*}
which is a morphism of Lie algebra crossed module.
}
\subsection{Inner automorphisms}
Recall the inner derivation Lie $2$-algebra of $\g$ defined in Section $3$, which is given by
\begin{equation*}
\CD
   \inn(\g): \inn^{-1}(\g)\triangleq \End^{-1}(\g)@>\bar{d}>>\inn^{0}(\g),
\endCD
\end{equation*}
where  $\inn^0(\g)=\Img \overline{\add}_0+ \Img\bar{d}$ is indeed the set of $1$-coboundaries of $\g$ with the coefficients in the adjoint representation. Moreover,
$\inn(\g)$ is an ideal of the derivation Lie $2$-algebra $\Der(\g)$.
The goal of this subsection is to define the inner automorphism $2$-group of a Lie $2$-algebra $\g$, which is a normal sub-crossed module \cite{Norrie} of $\Aut(\g)$. The main ingredient used in this subsection is the exponential map given in Corollary \ref{cor:e}.

For the Lie group $\Aut(\g)$, we have the adjoint map $\Ad: \Aut(\g)\longrightarrow \Aut(\Aut(\g))$ defined by
\begin{equation*}
\left\{\begin{array}{rcl} \Ad(A)(A')&=&A\cdot A'\cdot A^{-1}=A\diamond A'\diamond A^{-1},
\\ \Ad(A)(\tau)&=&A\cdot \tau\cdot A^{-1}=A\triangleright \tau,\\
\Ad(\tau)(A)&=&\tau\cdot A\cdot \tau^{-1}=(A,\tau\star(A\triangleright \tau^{-1})),\\
\Ad(\tau)(\tau')&=&\tau\cdot \tau'\cdot \tau^{-1}=\tau\star \tau'\star \tau^{-1}
.\end{array}\right.
\end{equation*}
We use the same notation to denote the induced adjoint representation of $\Aut(\g)$ on its Lie algebra $\Der(\g)$. Concretely, we have
\begin{lem}
The map $\Ad: \Aut(\g)\longrightarrow GL(\Der(\g))$ is formulated as follows:
\begin{equation*}
\left\{\begin{array}{rcl} \Ad(A)(X,l_X)&=&(A_0\circ X_0\circ A^{-1}_0, A_1\circ X_1\circ A^{-1}_1, l_{A\circ X\circ A^{-1}}),
\\ \Ad(A)(\Theta)&=&A_1\circ\Theta\circ A_0^{-1},\\
\Ad(\tau)(X,l_X)&=&((X,l_X),X_1\circ \tau^{-1}+\tau\circ X_0+\tau\circ d\circ X_0\circ \tau^{-1}),\\
\Ad(\tau)(\Theta)&=&(I+\tau\circ d)\circ \Theta \circ (I+d\circ \tau)^{-1},\end{array}\right.
\end{equation*}
where $l_{A\circ X\circ A^{-1}}$ is defined by
    \begin{eqnarray}\label{equaiii}l_{A\circ X\circ A^{-1}}(x,y)\nonumber&=&A_1l_X(A_0^{-1}x,A_0^{-1}y)-A_1X_1A_1^{-1}A_2(A_0^{-1}x,A_0^{-1}y)
    \\ &&+A_2(X_0A_0^{-1}x,A_0^{-1}y)+A_2(A_0^{-1}x,X_0A_0^{-1}y).
    \end{eqnarray}
\end{lem}
For any Lie group $K$ with Lie algebra $\mathfrak{k}$, we have
$$k\cdot e^{X}\cdot k^{-1}=e^{\Ad(k)X},\ \ \ \ \ \ \  \forall k\in K,X\in \mathfrak{k}.$$
Depending on this, we can derive several interesting facts.

\begin{cor}\label{most}
For all $(X,l_X)\in \Der^0(\g), \Theta\in \Der^{-1}(\g), A\in \Aut^0(\g)$ and $\tau\in \Aut^{-1}(\g)$, we have
\begin{enumerate}
 \item[\rm(i)]  $A\diamond e^{(X,l_X)}\diamond A^{-1}=e^{(A\circ X\circ A^{-1},l_{A\circ X\circ A^{-1}})}$, where $l_{A\circ X\circ A^{-1}}$ is given by \eqref{equaiii};
\item[\rm(ii)] $\tau\star e^\Theta\star \tau^{-1}=e^{(I+\tau\circ d)\circ \Theta \circ (I+d\circ \tau)^{-1}};$
\item[\rm(iii)] $A\triangleright e^\Theta=e^{A_1\circ\Theta\circ A_0^{-1}};$
\item[\rm(iv)] $\tau\star (e^{(X,l_X)}\triangleright \tau^{-1})=e^{X_1\circ \tau^{-1}+\tau\circ X_0+\tau\circ d\circ X_0\circ \tau^{-1}}.$
\end{enumerate}
\end{cor}
\begin{cor}
For all $x\in \g_0, \Theta\in \Der^{-1}(\g)$ and $A\in \Aut^0(\g)$, we have
\begin{eqnarray}
\label{equaa}A\diamond e^{\bar{d}(\Theta)}\diamond A^{-1}&=&e^{\bar{d}(A_1\circ\Theta\circ A_0^{-1})},\\
\label{equab}A\diamond e^{\overline{\add}_0(x)}\diamond A^{-1}&=&e^{\overline{\add}_0(A_0x)+\bar{d}(A_2(x,A_0^{-1}\cdot))}.
\end{eqnarray}
\end{cor}
\pf Note that $\bar{d}(\Theta)=(\delta(\Theta), l_{\delta (\Theta)})$. Firstly, it is simple to check that
\begin{eqnarray*}
A\circ \delta(\Theta)\circ A^{-1}=\delta(A_1\circ \Theta \circ A_0^{-1}).
\end{eqnarray*}
Then, by (\ref{equaiii}) and the fact that $A,A^{-1}$ are Lie $2$-algebra homomorphisms, we find that
\begin{eqnarray*}
&&l_{\delta(A_1\circ\Theta\circ A_0^{-1})}(y,z)\\ &=&A_1\Theta A_0^{-1}[y,z]-[A_1\Theta A_0^{-1}y,z]-[y,A_1\Theta A_0^{-1}z]
\\&=&A_1\Theta[A_0^{-1}y,A_0^{-1}z]-A_1\Theta dA_1^{-1}A_2(A_0^{-1}y,A_0^{-1}z)-A_1[\Theta A_0^{-1}y,A_0^{-1}z]+A_2(d\Theta A_0^{-1}y, A_0^{-1}z)\\ &&-A_1[A_0^{-1}y,\Theta A_0^{-1}z]+A_2(A_0^{-1}y,d\Theta A_0^{-1}z)\\ &=&A_1l_{\delta(\Theta)}(A_0^{-1}y, A_0^{-1}z)-A_1\Theta dA_1^{-1}A_2(A_0^{-1}y,A_0^{-1}z)+A_2(d\Theta A_0^{-1}y, A_0^{-1}z)+A_2(A_0^{-1}y,d\Theta A_0^{-1}z)\\ &=&l_{A\circ \delta(\Theta)\circ A^{-1}}(y,z).
\end{eqnarray*}
Hence, by (i) of Corollary \ref{most}, we obtain (\ref{equaa}).

Note that $\overline{\add}_0(x)=(\ad_0(x),l_3(x,\cdot,\cdot)).$ Since $A$ is a Lie $2$-algebra homomorphism, it is simple to verify
$$A\circ \ad_0(x)\circ A^{-1}=\ad_0(A_0x)+\delta(A_2(x,A_0^{-1}\cdot)).$$
Set $\Lambda:=A_2(x,A_0^{-1}\cdot)$ and $\ad_x:=\ad_0(x)$ for simplicity. Using again the fact that $A$ and $A^{-1}$ are Lie $2$-algebra homomorphisms, we have
\begin{eqnarray*}
(l_3(A_0x,\cdot,\cdot)+l_{\delta(\Lambda)})(y,z) &=&l_3(A_0x,y,z)+A_2(x,[A_0^{-1}y,A_0^{-1}z]-dA_1^{-1}A_2(A_0^{-1}y,A_0^{-1}z))\\ &&-[A_2(x,A_0^{-1}y),z]-[y,A_2(x,A_0^{-1}z)]\\ &=&A_2(A_0^{-1}y,\ad_xA_0^{-1}z)+A_2(\ad_xA_0^{-1}y,A_0^{-1}z)-\ad_{A_0x}A_2(A_0^{-1}y,A_0^{-1}z)\\ &&-A_2(x,dA_1^{-1}A_2(A_0^{-1}y,A_0^{-1}z))+
A_1l_3(x,A_0^{-1}y,A_0^{-1}z)\\ &=&-(\ad_{A_0x}+\Lambda d)A_2(A_0^{-1}y,A_0^{-1}z)+A_1l_3(x,A_0^{-1}y,A_0^{-1}z)\\ &&+A_2(A_0^{-1}y,\ad_xA_0^{-1}z)+A_2(\ad_xA_0^{-1}y,A_0^{-1}z)\\ &=&l_{A\circ \overline{\add}_0(x)\circ A^{-1}}(y,z),
\end{eqnarray*}
where the second equality follows from the condition (iii) of Definition \ref{defi:Lie 2 homo}.
 By (i) of Corollary \ref{most}, we obtain (\ref{equab}).\qed\vspace{3mm}

Denote by $\Inn^0(\g)$ and $\Inn^{-1}(\g)$ the connected Lie groups $e^{\inn^0(\g)}$ and $ e^{\inn^{-1}(\g)} $ respectively. Then we get a natural sub-crossed module of $\Aut(\g)$:
\begin{equation*}
\CD
   \Inn(\g): \Inn^{-1}(\g) @>\partial>> \Inn^0(\g).
\endCD
\end{equation*}
\begin{thm}
With the notations above,
$\Inn(\g)$ is a normal sub-crossed module of $\Aut(\g)$.
\end{thm}
Since $\Inn(\g)$ is the image of the inner derivation Lie 2-algebra under the exponential map, we call it the {\bf inner automorphism $2$-group} of $\g$.

For a strict Lie 2-algebra $\g=(\g_0\oplus\g_{-1},d,[\cdot,\cdot])$, consider the sub-complex $\Sinn(\g)$ of $\SDer(\g)$ given by $\Sinn^0(\g)=\Img(\ad_0)$ and  $\Sinn^{-1}(\g)=\Img(\ad_1)$. It is straightforward to see that $\Sinn(\g)$ is an ideal of $\SDer(\g)$. Denote by $\SInn^0(\g)$ and $\SInn^{-1}(\g)$ the connected Lie groups $e^{\Sinn^0(\g)}$ and $ e^{\Sinn^{-1}(\g)} $ respectively. Then $\SInn(\g)$ is a normal sub-crossed module of $\SAut(\g)$, which is consistent with the inner automorphism $2$-group given in \cite{Norrie}.
\begin{rmk}{\rm
  In \cite{Urs}, an inner automorphism $3$-group of a strict Lie $2$-group is given from a slight different viewpoint. Their inner automorphism $3$-group can remember the center, which is different from ours.
 \emptycomment{
\left\{\begin{array}{rcl} Obj(\mathbb{G})&=&G_0,
\\ Mor_1(\mathbb{G})&=&\{(\alpha,\xi,\beta)|\alpha,\beta\in G_0,\xi\in G_{-1}\},\\
Mor_2(\mathbb{G})&=&G_{-1},\end{array}\right.
\end{equation*}
where a $1$-morphism $\beta\longrightarrow \gamma$ is exactly
$$(\alpha,\xi,\beta): \beta\longrightarrow d^{\mathbb G}(\xi)\alpha\beta,$$
 and a $2$-morphism $(\alpha,\xi,\beta)\longrightarrow (\alpha',\xi',\beta)$ is
 $$\eta:(\alpha,\xi,\beta)\longrightarrow (\alpha',\xi\eta^{-1},\beta).$$
 }
  }
\end{rmk}

 \end{document}